\documentclass[12pt]{article}
\usepackage[british]{babel}
\usepackage{authblk}
\usepackage{amsmath,hyperref}
\usepackage{amssymb}
\usepackage{enumerate}
\usepackage[amsmath,thmmarks]{ntheorem}

\newtheorem{proposition}{Proposition}

\usepackage{microtype}
\newtheorem{rem}{Remark}

\usepackage{setspace}
\usepackage{subcaption}
\usepackage{multicol}
\usepackage{booktabs}
\usepackage{multirow}
\usepackage[dvipsnames]{xcolor}
\usepackage[compact]{titlesec}
\usepackage{graphicx} 
\usepackage{amsmath}
\usepackage{amssymb}
\usepackage{algorithm}
\usepackage{algorithmic}
\usepackage{mathtools}
\usepackage{graphicx} 
\usepackage{listofitems} 
\usepackage{tikz}
\usetikzlibrary{shapes.geometric, arrows}
\usepackage{listofitems} 
\usetikzlibrary{arrows.meta} 
\usepackage[outline]{contour} 
\usepackage{xcolor}
\colorlet{myred}{red!80!black}
\colorlet{myblue}{blue!80!black}
\colorlet{mygreen}{green!60!black}
\colorlet{myorange}{orange!70!red!60!black}
\colorlet{mydarkred}{red!30!black}
\colorlet{mydarkblue}{blue!40!black}
\colorlet{mydarkgreen}{green!30!black}

\let\-=\mathbf
\tikzset{
  >=latex, 
  node/.style={thick,circle,draw=myblue,minimum size=22,inner sep=0.5,outer sep=0.6},
  node in/.style={node,green!20!black,draw=mygreen!30!black,fill=mygreen!25},
  node hidden/.style={node,blue!20!black,draw=myblue!30!black,fill=myblue!20},
  node convol/.style={node,orange!20!black,draw=myorange!30!black,fill=myorange!20},
  node out/.style={node,red!20!black,draw=myred!30!black,fill=myred!20},
  connect/.style={thick,mydarkblue}, 
  connect arrow/.style={-{Latex[length=4,width=3.5]},thick,mydarkblue,shorten <=0.5,shorten >=1},
  node 1/.style={node in}, 
  node 2/.style={node hidden},
  node 3/.style={node out}
}

\theoremstyle{nonumberplain}
\theoremheaderfont{\normalfont\itshape}
\theorembodyfont{\normalfont}
\theoremseparator{.}
\theoremsymbol{\ensuremath{\square}}

\renewcommand{\Gamma}{\varGamma}
\renewcommand{\epsilon}{\varepsilon}

\renewcommand{\leq}{\leqslant}

\newcommand{\bdmu}{\boldsymbol{\mu}}
\newcommand{\bdtheta}{\boldsymbol{\theta}}

\newcommand{\bmu}{\boldsymbol{\mu}}
\newcommand{\bpsi}{\boldsymbol{\psi}}

\title{Nemytskii neural operator: a nonlinear model reduction method for parametrized partial differential equations}

\author[1]{Jingye Li}

\author[1]{Alex Bespalov}

\author[1]{Jinglai Li}

\affil[1]{School of Mathematics, University of Birmingham, Birmingham, United Kingdom}

\date{}
\begin{document}
\maketitle

\begin{abstract} 
We introduce a Nemytskii neural operator framework for nonlinear model reduction of parametrized steady-state partial differential equations. The method generalizes reduced basis approaches by replacing linear combinations of basis functions with a structured nonlinear mapping realized through a pointwise Nemytskii operator acting on fixed feature functions. Feature functions are learned offline via nonlinear dimension reduction from high-fidelity snapshots, and a hypernetwork maps model parameters to a lightweight reconstruction network, which is further refined online using physics-informed residual minimization. The Nemytskii structure preserves analytical regularity and enables efficient evaluation of spatial and parametric derivatives, 
leading to fast online adaptation. 
Numerical experiments demonstrate that the proposed method consistently outperforms linear model reduction techniques, particularly for complex solution manifolds.
\end{abstract}


\section{Introduction}
\numberwithin{equation}{section}

\label{cha:introduction}

Parametrized partial differential equations (PDEs) arise in a wide range of scientific and engineering applications, including uncertainty quantification, optimization, and control, where one must repeatedly solve the same PDE for many different parameter values \cite{hesthaven2016certified,quarteroni2015reduced}.  High-fidelity discretizations such as finite element, spectral, or finite difference methods can be prohibitively expensive in this many-query context.  
The reduced basis (RB) methods \cite{hesthaven2016certified,quarteroni2015reduced} 
 address this challenge by assuming that the solution manifold can be well approximated by a low-dimensional linear subspace of the solution space.
In an expensive offline stage, RB methods construct a reduced approximation space (typically from high-fidelity snapshots); in the online stage, they compute solutions by solving the governing equations projected onto this reduced space. Because the online solves involve only a small number of degrees of freedom,  they are often significantly more efficient than solving the full PDE directly. 

While the standard RB methods can yield significant speedups, they still require the solution of algebraic systems at the online stage.
This step can remain computationally demanding, particularly when the underlying PDE is nonlinear.
Non‐intrusive RB approaches reduce computational cost during the online stage by replacing the traditional solver with a regression model that is constructed in the offline stage and directly maps parameters to RB coefficients. While noting that early works employ conventional interpolation methods (see, e.g., \cite{casenave2015nonintrusive,xiao2016non}), 
our focus here is on more recent developments where trained neural networks (NN) are utilized to predict RB coefficients from parameters (see, e.g., \cite{hesthaven2018non,wang2019non,chen2024gpt}).

It is important to emphasize that all RB-type methods rely on the central assumption that the solution manifold admits a good approximation by a low-dimensional linear subspace; in practice, this assumption may not hold, or the dimension of the required linear subspace may be unacceptably high.
When linear approximation is not effective,  nonlinear reduced-order modeling (ROM) provides a natural alternative framework.
These methods (e.g.,  \cite {lee2020model,fresca2021comprehensive}) seek to represent the solution set using nonlinear low-dimensional structures, which can capture complex solution behavior more efficiently than linear subspaces.
In contrast to RB methods, which construct reduced spaces from explicitly defined basis functions and thus impose a clear mathematical structure, many existing nonlinear ROM approaches rely on fully data-driven models in which the solution representation and its structure are learned primarily from data.
As a result, such methods often require large training datasets and complex neural networks, leading to increased training cost and limited robustness and generalization.
A more detailed discussion of nonlinear ROM methods and their limitations is provided in Section~\ref{sec:rom}.

Motivated by these considerations, we introduce a Nemytskii neural operator (NNO) framework for nonlinear ROM of parametrized PDEs. 
The proposed approach can be regarded as a nonlinear generalization of RB methods:
instead of representing the solution as a linear combination of basis functions, it constructs the approximation through a nonlinear mapping realized by a pointwise Nemytskii  operator acting on fixed feature functions.
The pointwise (superposition) nature of the Nemytskii operator ensures that the solution is constructed locally from the fixed feature functions; as a result, the resulting nonlinear mapping is highly structured and principled, with its structure strongly constrained by the feature functions.
This property is closely analogous to the role played by basis functions in  RB methods.
As such the proposed NNO framework overcomes the limitations of linear approximation in RB methods, 
while retaining their ability to impose explicit analytical structure through fixed feature functions, a key distinction from generic nonlinear ROM.

Beyond the Nemytskii representation, the proposed NNO framework consists of two additional key components that are essential to its performance: online adaptation and NNO-based dimension reduction.

Many learning–based ROM/RB methods, such as non-intrusive RB approaches, rely on surrogate models that are constructed entirely during the offline stage. In the online stage, the parameter value of interest is simply provided as input to the trained model to obtain an approximate solution.
The surrogate models constructed offline may not be sufficiently accurate, particularly in parameter regimes that are poorly represented in the training data,
leading to unreliable solution predictions with no correction mechanism available during the online stage. 
To address this issue, the proposed NNO framework incorporates an online adaptation mechanism that allows the solution approximation to be adapted for each new parameter instance during the online stage.
The adaptation is performed in an unsupervised manner using physics-based residuals associated with the governing PDE \cite{raissi2019physics,cai2021physics}.
Importantly the NNO framework allows us to represent the solution using a lightweight network defined on the feature space,
enabling efficient adaptation during the online stage.

Another important component of the proposed method is nonlinear dimension reduction (feature extraction) performed at the offline stage. The purpose of this step is to learn a set of  feature functions from high-fidelity solution snapshots.
In the proposed method, the feature functions are also constructed through a Nemytskii-operator–based nonlinear mapping applied to the snapshot data.
This design enables the model to effectively learn the underlying analytical structure from the solution snapshots.
The learned feature functions are fixed during the online stage to ensure computational efficiency.

In summary, the main contribution of this paper is the introduction of a Nemytskii neural operator (NNO) framework for constructing a nonlinear approximation of the solution manifold of parametrized PDEs. The proposed approach combines structured nonlinear solution representation with offline nonlinear dimension reduction and efficient online adaptation, enabling expressive approximation while retaining analytical structure and computational efficiency. By exploiting the pointwise nature of the Nemytskii operator and operating in a low-dimensional feature space, the method provides a principled alternative to both linear RB methods and generic data-driven nonlinear ROM approaches.

The remainder of this work is organized as follows. Section \ref{sec:background} presents background and preliminary material, including the parametrized PDE model of interest, standard and non-intrusive RB methods, nonlinear ROM and
physics-informed neural networks (PINNs). Section~\ref{sec:method} introduces the proposed NNO method, outlining its offline and online components in detail. Section \ref{sec:examples} presents numerical experiments that demonstrate the effectiveness of the proposed method. Finally, Section \ref{sec:conclusion} concludes the paper with a summary and potential directions for future research.

\section{Preliminaries and background} \label{sec:background}
\numberwithin{equation}{section}
\subsection{Problem Setup} \label{sec:setup}
We consider a generic parametrized steady-state PDE system:
\begin{subequations}\label{eq:non_pde}
\begin{equation}\label{eq:mainpde}
    {\mathcal{N}}\bigl(u(\mathbf{x},\boldsymbol{\mu}),\boldsymbol{\mu}\bigr)=0, \quad \mathbf{x} \in \Omega,
\end{equation}
subject to Dirichlet (essential) boundary condition (BC):
\begin{equation}\label{e:bc}
u(\-x,\bmu)=b(\-x,\bmu)\quad \mathbf{x} \in \partial\Omega.
\end{equation}
\end{subequations}
 Here, ${\mathcal{N}}$ is a (possibly nonlinear) differential operator, $u(\mathbf{x},\boldsymbol{\mu})$ is the solution, 
 $\-x$ is the $d$-dimensional spatial variable,
 and $\Omega \subset \mathbb{R}^d$ is the spatial domain with $\partial\Omega$ being its boundary; 
 $\boldsymbol{\mu} \in {Y}$ represents the collection of model parameters, with ${Y}$ being the parameter domain. 
 Since natural boundary conditions can be incorporated implicitly in the weak formulation and therefore do not need to be treated explicitly in the methods considered in this work, we restrict our attention to Dirichlet BCs.
 
 In many practical settings, $\boldsymbol{\mu}$ is modelled as a random variable with a prescribed probability distribution $\pi_{\boldsymbol{\mu}}(\cdot)$. 
 Our objective is to numerically solve the parametrized PDE for a given value of $\boldsymbol{\mu}$.
Conventionally, such PDEs are solved using standard numerical schemes such as the finite difference method (FDM), the finite element method (FEM), or spectral methods. Here, we present the Galerkin FEM as an example (see, e.g., \cite{hughes2003finite}) .
Let $\mathcal V$ be a Hilbert space equipped with the inner product 
$\langle\cdot,\cdot\rangle_{\mathcal V}$ and norm $\|\cdot\|_{\mathcal V}$. 
We assume that $\mathcal V$ is chosen to satisfy the Dirichlet boundary conditions.
Its dual space is denoted by $\mathcal{V}'$, and the duality pairing is given by $\langle \cdot,\cdot \rangle_{\mathcal{V}',\mathcal{V}} : \mathcal{V}' \times \mathcal{V} \rightarrow \mathbb{R}$. We can view the differential operator ${\mathcal{N}}$ in \eqref{eq:non_pde} as a mapping:
\begin{equation*}
{\mathcal{N}}: \mathcal{V} \times {Y} \longrightarrow \mathcal{V}'.
\end{equation*}
The weak formulation of \eqref{eq:non_pde} is: 
\begin{equation}\label{non_pde_sol2}
\langle {{\mathcal{N}}}\bigl(u(\cdot,{\boldsymbol{\mu}}),{\boldsymbol{\mu}}\bigr), v \rangle_{\mathcal{V}',\mathcal{V}} = 0, \quad \forall v \in \mathcal{V},
\end{equation}
and a solution to Eq.~\eqref{non_pde_sol2} is called a weak solution to Eq.~\eqref{eq:non_pde}.

Next the weak solution $u \in \mathcal{V}$ is constructed in a finite-dimensional subspace $\mathcal{V}_n \subset \mathcal{V}$, spanned by $n$ basis functions $\{\eta_1,\eta_2,\ldots,\eta_n\}$. The numerical solution is expressed as a linear combination:
\begin{equation}\label{high_fide_sol_fix}
u_h(\-x,{\boldsymbol{\mu}}) = \sum_{j=1}^{n} \tilde{u}_j({\boldsymbol{\mu}})\,\eta_j(\-x),
\end{equation}
and we seek the coefficient vector
\begin{equation}\label{e:sol_coef_full}
\tilde{\mathbf{u}}({\boldsymbol{\mu}}) = \begin{bmatrix} \tilde{u}_1({\boldsymbol{\mu}}) & \ldots & \tilde{u}_n({\boldsymbol{\mu}}) \end{bmatrix}^T \in \mathbb{R}^{n},
\end{equation}
such that
\begin{equation}\label{non pde sol fem 1}
\langle{\mathcal{N}}\bigl( u_h(\cdot,{\boldsymbol{\mu}}),{\boldsymbol{\mu}}\bigr),\eta_j \rangle_{\mathcal{V}',\mathcal{V}}=0, \,j = 1,\ldots,n.
\end{equation}
Next we need to solve the equation system \eqref{non pde sol fem 1}. If ${\mathcal{N}}$ is a linear operator, it reduces to a linear system. If ${\mathcal{N}}$ is nonlinear, 
a common approach is to solve the equation system using a linearization procedure (e.g., Newton's method \cite{demmel1997applied}). 

The presumption here is that we use a sufficiently large number of basis functions so that the numerical approximate solution $u_h$ is of high accuracy
and is referred to as the \emph{high-fidelity solution} in this paper.  
To achieve such accuracy, the number of basis functions required is often very large, which in turn leads to a very 
high-dimensional algebraic equation system in Eq.~\eqref{non pde sol fem 1}. 
Assembling and solving such a large system---especially using iterative linearization in the case of a nonlinear operator---becomes computationally expensive in terms of both memory and CPU time.
To this end, the RB method is a popular approach for solving parametrized PDEs, as it can often significantly reduce
the dimensionality of the algebraic equation system. 

\subsection{Standard reduced basis method}
\label{cha:RB sec}

In this section, we give a brief overview of the RB method. 
The main idea of the method is rather straightforward: first, we construct a low-dimensional approximation space that can well approximate the solution set across the parameter space, and then we solve the PDE for a specific parameter value in this reduced low-dimensional space, significantly lowering computational costs compared to the standard approaches. 
The standard RB method is implemented in two stages: an \emph{offline stage} and an \emph{online stage}.
In the offline stage, a reduced-dimensional space for solution approximation is constructed. This stage is computationally expensive but performed only once. In the online stage, the reduced space obtained offline is used to rapidly approximate solutions for new parameter values,  enabling real-time simulations of the PDE model.
\subsubsection{Low-dimensional linear subspace}
Let  the solution set of Eq.~\eqref{non_pde_sol2} be, 
\[
\mathcal M
:=
\{\, u(\cdot,\bmu) \in \mathcal V \mid \bmu \in  Y \,\} \subset \mathcal{V}, 
\]
and we assume that $\mathcal M$ can be well approximated 
by a $r$-dimensional linear subspace:
\begin{equation}\label{rb_space}
\mathcal{V}_{r} = \mathrm{span}\{\psi_1, \psi_2, \dots, \psi_r\},
\end{equation}
where $\psi_1,...,\psi_r$ are the basis functions spanning $\mathcal V_r$.
That is, each $u(\cdot,\bmu) \in \mathcal M$ can be well approximated by a linear combination 
of $\psi_1,...,\psi_r$. For RB methods to be effective, the reduced dimension $r$ needs to be much smaller 
than $n$, the dimensionality of the FEM space.
A key step in RB methods is to obtain the basis functions $\psi_1,...,\psi_r$,
which are constructed in the offline stage. 

\subsubsection{Offline stage}
The first step in the offline stage is to compute  a set of high-fidelity solutions, often called snapshots,  corresponding to some selected parameter values. 
In practice, one usually randomly generates the parameter values from the underlying distribution of the parameter: 
\[\{\boldsymbol{\mu}^{(1)}, \boldsymbol{\mu}^{(2)}, \dots, \boldsymbol{\mu}^{(N_s)}\},\,\,\mathrm{with}\,\,
\boldsymbol{\mu}^{(i)} \sim \pi_{\boldsymbol{\mu}}(\cdot)\,\,\mathrm{for}\,\, i=1,\ldots,N_s.\]
For each parameter value $\boldsymbol{\mu}^{(i)}$, we compute the high-fidelity solution using a standard numerical method (e.g., the Galerkin method discussed in Section \ref{sec:setup}). 
Let the high-fidelity solutions for each parameter value be denoted as
\begin{equation}\label{discrete_sol}
{u}_h(\cdot,\boldsymbol{\mu}^{(i)}), \quad i=1,2,\dots,N_s;
\end{equation}
these are represented as in Eq.~\eqref{high_fide_sol_fix}, and the coefficient vectors 
$\tilde{\mathbf{u}}_h(\boldsymbol{\mu}^{(i)})$
are given by Eq.~\eqref{e:sol_coef_full}.
We reinstate that $\tilde{\mathbf{u}}_h(\boldsymbol{\mu}^{(i)})$ are $n$-dimensional vectors, and 
they form the snapshot matrix $\mathbf{S}$:
\begin{equation}\label{snapshot_matrix}
\mathbf{S} = [\tilde{\mathbf{u}}_h(\boldsymbol{\mu}^{(1)}) \quad \tilde{\mathbf{u}}_h(\boldsymbol{\mu}^{(2)}) \quad \dots \quad 
\tilde{\mathbf{u}}_h(\boldsymbol{\mu}^{(N_s)})] \in \mathbb{R}^{n \times N_s}.
\end{equation}
As an alternative to random sampling, greedy algorithms are also often
employed to adaptively generate snapshots; see, e.g.,
\cite{quarteroni2015reduced,hesthaven2016certified}.

Next, we build the reduced space using these snapshots, often via the {proper orthogonal decomposition} (POD) technique \cite{hesthaven2016certified,quarteroni2015reduced}. 
Specifically, we perform the singular value decomposition (SVD) of the snapshot matrix $\mathbf{S}$:
\begin{equation}\label{svd}
\mathbf{S} = \mathbf{B} \mathbf{\Sigma} \mathbf{V}^T,
\end{equation}
where $\mathbf{B} \in \mathbb{R}^{n \times n}$, $\mathbf{\Sigma}\in \mathbb{R}^{n\times N_s}$ (diagonal matrix containing singular values), and $\mathbf{V} \in \mathbb{R}^{N_s \times N_s}$.
We select the first $r \ll n$ singular vectors corresponding to the largest singular values to construct the basis functions of $\mathcal V_r$:
\begin{equation}\label{RB_basis}
\psi_i(\mathbf{x}) = \sum_{j=1}^{n} \mathbf{B}_{j,i}\,\eta_j(\mathbf{x}), \quad i=1,\dots,r,
\end{equation}
where $\{\eta_j(\mathbf{x})\}_{j=1}^{n}$ are the original basis functions from the high-fidelity solution space. 


\subsubsection{Online stage}\label{sec:online_stage}

Once the reduced space $\mathcal{V}_r = \text{span}\{\psi_1, \dots, \psi_r\}$ has been constructed in the offline stage, the online stage aims to compute an efficient approximation of the high-fidelity solution for a new parameter value ${\boldsymbol{\mu}}^* \in {Y}$.
We seek a reduced solution $u_{\rm rb}(\-x, \boldsymbol{\mu}^*) \in \mathcal{V}_r$ expressed as a linear combination of the reduced basis functions:
\begin{equation}\label{e:rb_sol}
   u_{\rm rb}(\mathbf{x},\boldsymbol{\mu}^*) = \tilde{{u}}_{1}({\boldsymbol{\mu}}^*)\psi_1(\mathbf{x}) + \ldots + \tilde{{u}}_{r}({\boldsymbol{\mu}}^*)\psi_{r}(\mathbf{x}),
\end{equation}
and our goal now is to compute the coefficient vector:
\begin{equation}\label{e:sol_coef}
\tilde{\mathbf{u}}_{\rm rb}({\boldsymbol{\mu}}^*) = \begin{bmatrix} \tilde{u}_1({\boldsymbol{\mu}}^*) & \ldots & \tilde{u}_r({\boldsymbol{\mu}}^*) \end{bmatrix}^T \in \mathbb{R}^{r}.
\end{equation}
Here we can use, for example, the Galerkin projection method, except that we now construct the solution in the reduced space $\mathcal{V}_r$,
rather than in $\mathcal{V}_n$. Namely, Eq.~\eqref{non pde sol fem 1} becomes
\begin{equation}\label{e:projection_rb}
\langle{\mathcal{N}}\bigl( u_{\rm rb}(\cdot,{\boldsymbol{\mu}}^*),{\boldsymbol{\mu}}^*\bigr),\psi_j \rangle_{\mathcal{V}',\mathcal{V}}=0,\quad j = 1,\ldots,r,
\end{equation}
where $u_{\rm rb}$ is given by Eq.~\eqref{e:rb_sol}. 
Therefore, using the reduced bases, we obtain a 
system of $r$ algebraic equations, compared to a system of  $n$ equations in the original Galerkin projection method.  In many practical problems, $r$ is much smaller than $n$, which can result in significant computational savings. When the underlying algebraic equations are linear and certain conditions are satisfied, the resulting linear system can be solved very efficiently, despite the fact that the associated system matrix is, in general, dense. Even in nonlinear cases, the method still improves efficiency due to the substantially reduced system size (see, e.g., \cite{carlberg2011efficient,chaturantabut2010nonlinear}).

\subsection{Non-intrusive RB methods}
In the standard RB method, especially when the system is nonlinear, equations~\eqref{e:projection_rb} need to be solved numerically at the online stage. 
This reduced system is much less expensive than the original system, but may still require rather substantial computational time. 

To this end, the non-intrusive RB methods are proposed to further reduce the online computational cost. 
The basic idea behind these methods is rather straightforward: instead of solving for the coefficient at the online stage, 
one imposes a regression model representing the coefficients $\tilde{\mathbf{u}}_{\rm rb}$ as a function of the model parameter $\boldsymbol{\mu}$,
which is trained at the offline stage with additional high-fidelity solution data. 
We refer to \cite{yu2019non,padula2024brief} for more comprehensive reviews of the non-intrusive RB methods. 
 The non-intrusive RB methods do not require solving the {reduced problem} in the {online} stage,
 and thus make the online computation much faster. 
As mentioned earlier, we here focus on methods that use NN models to compute the reduced solutions,
while noting that other surrogate models are also available \cite{casenave2015nonintrusive,guo2018reduced,guo2019data}. In \cite{hesthaven2018non,chen2018greedy}, high-fidelity solutions for a set of parameter values are used as labeled data to train an NN with parameter as input and coefficient vector as output.
The basic procedure of the method is the following: 
  \begin{enumerate}
 \item Generate the reduced basis $\{\psi_1,\ldots,\psi_r\}$ using the same procedure as in the standard RB method.
 \item Write the RB solution as Eq.~\eqref{e:rb_sol} with the coefficient $\tilde{\mathbf{u}}_{\rm rb}$ given by Eq.~\eqref{e:sol_coef}.
 \item Assume an NN model  $\tilde{\mathbf{u}}_{\rm rb}(\boldsymbol{\mu}) = f(\boldsymbol{\mu};\boldsymbol{\theta}_f)$, where $\boldsymbol{\theta}_f$ are the model parameters.
 \item Sample a new set of  $N$ parameter values $\boldsymbol{\mu}^{(1)}, \ldots, \boldsymbol{\mu}^{(N)}$ and compute the associated
 high-fidelity solutions $u_h(\-x,\boldsymbol{\mu}^{(1)}), \ldots, u_h(\-x,\boldsymbol{\mu}^{(N)})$.
\item Construct the labeled training data: \ $(\boldsymbol{\mu}^{(1)},
\tilde{\-u}^{(1)}),\ldots,
(\boldsymbol{\mu}^{(N)},\tilde{\-u}^{(N)})$,
 where $\tilde{\-u}^{(i)}$ denotes the coefficients of 
 $u_h(\-x,\boldsymbol{\mu}^{(i)})$ projected onto the reduced basis functions $\psi_1, \dots, \psi_r$.
 \item Estimate $\boldsymbol{\theta}_f$ by minimizing the supervised loss function:
 \[ L_{\text{RB}}(\boldsymbol{\theta}_f) = \sum_{i=1}^{N} \|\tilde{\-u}^{(i)}-f(\boldsymbol{\mu}^{(i)};\boldsymbol{\theta}_f)\|_2^2.\]
 \end{enumerate}
 Several variants of the method have been developed as well. 
 For example, the number of labelled data required to train the network can be reduced if the physics-informed loss is added to the training loss \cite{chen2021physics}.
Other variants of NN based non-intrusive RB include
 \cite{fresca2022pod}, which constructs the NN model using autoencoder architectures,
 \cite{pichi2024graph}, where graph-based approaches have been proposed to handle unstructured spatial domains.

\subsection{Nonlinear reduced-order modeling}\label{sec:rom}
A major limitation of RB methods is their reliance on linear approximation subspaces, which can require a prohibitively large number of basis functions for solution manifolds with slowly decaying Kolmogorov $n$-widths.
Nonlinear reduced-order modeling (ROM) addresses this limitation by approximating the solution set using low-dimensional nonlinear manifolds rather than linear subspaces.
Many such nonlinear ROM techniques can be 
interpreted as a generative model.
Namely, we assume that the solution set $\mathcal M$ can be well approximated by a low-dimensional  manifold
\begin{equation} \mathcal{M}_G := \{ {G}(\-z) \mid \-z \in Z \subset\mathbb{R}^\ell \}, \label{e:MG}
\end{equation}
where $\-z\in Z$ is a $\ell$-dimensional latent variable, and 
$G:$ $Z\rightarrow \mathcal{V}$
is a mapping from the latent space $Z$ to function space $\mathcal V$. 
In practice, the mapping is often parametrized by a neural network
and then trained over a collection of high-fidelity solution snapshots \cite{fresca2021comprehensive}. 

Despite their increased expressive power, many nonlinear ROM approaches exhibit several notable
limitations.
First, they are typically trained in a purely data-driven manner and do not explicitly incorporate
the structure of the governing PDE and its solutions, which may lead to limited robustness and poor generalization
outside the training regime.
Second, because the approximation is learned entirely from data, nonlinear ROM methods typically
require substantially larger training data sets in order to accurately capture the geometry of the
solution manifold.
Finally, due to the high-dimensional nature of the mapping from the latent space to the solution
space, particularly on the output side, the generative networks employed in nonlinear ROM are often
large-scale, making training computationally expensive and memory-intensive.

\subsection{Physics-informed neural networks}
\label{cha:PINN}
In this section we review the physics-informed neural networks~(PINNs), another important ingredient of our method. 
PINNs~\cite{raissi2019physics} are a class of machine learning based methods for the numerical solution of PDEs.
Suppose we want to solve the boundary value problem ~\eqref{eq:non_pde}.
When using PINN, the solution to \eqref{eq:non_pde} is approximated by an NN $u_{\rm NN}(\mathbf{x};\boldsymbol{\theta})$
with $\boldsymbol{\theta}$ being the trainable network parameters. 
Since $u_{\rm NN}(\mathbf{x};\boldsymbol{\theta})$ is only an approximate solution of Eq.~\eqref{eq:non_pde},
it does not satisfy the PDE and boundary conditions exactly.
Instead, the residuals of the governing equations are used to guide training.

In this work, we adopt a weak-form (variational) formulation of PINNs.
Let $\{ v_i \}_{i=1}^k \subset \mathcal V$ denote a collection of $k$ selected test functions.
The weak residual associated with the PDE is defined as
\[
R(\bdtheta;v_i) = \left\langle 
{\mathcal N}\bigl(u_{\mathrm{NN}}(\cdot;\boldsymbol{\theta}),\boldsymbol{\mu}\bigr),\,
v_i
\right\rangle_{\mathcal V',\mathcal V}, \quad i=1,...,k.
\]
The basic idea here is that we should train the network $u_{\rm NN}(\mathbf{x};\boldsymbol{\theta})$ to minimize the residuals,
yielding loss function:
\begin{equation}\label{e:loss_pinn}
L_\mathrm{PINN}(\boldsymbol{\theta}) = \sum_{i=1}^k R^2(\bdtheta;v_i)
+ 
  \lambda_{BC} L_{BC}(\bdtheta)
    \end{equation}
where 
$\lambda_{BC}$ is a penalty coefficient, and $L_{BC}(\cdot)$ denotes the boundary-condition penalty:
\[
L_{BC}(\bdtheta)=  \| 
u_{\mathrm{NN}}(\mathbf{x};\boldsymbol{\theta}) - b(\mathbf{x},\bmu)
\|^2_{L^2(\partial\Omega)}.
\]
In practice, both terms in the loss function are approximated numerically using quadrature rules, leading to a discrete loss function that can be efficiently evaluated and minimized using stochastic gradient–based optimization methods.


To approximate the solution, PINNs typically employ a Multi-Layer Perceptron (MLP) network, consisting of an input layer, multiple hidden layers, and a single output layer~\cite{goodfellow2016deep}. Each layer is made up of a number of neurons which are connected to those in adjacent layers. 
PINNs often employ deep architectures with a substantial number of neurons to effectively capture complex solution structure. 
For example, in their seminal work  \cite{raissi2019physics}, the authors utilized a neural network with 9 layers, each containing 20 neurons, to solve the Navier–Stokes equation, and one with 4 hidden layers of 200 neurons each for the Allen–Cahn equation.

\section{Nemytskii neural operator with online adaptation} \label{sec:method}

\numberwithin{equation}{section}

\label{cha:Model sec}

\subsection{Nemytskii neural operator}

Similar to the nonlinear ROM, our method also aims to construct a nonlinear solution manifold,
that approximates the solution set of the parametrized PDE \eqref{eq:non_pde}.
We assume that we have a collection of ``feature functions'' $\psi_1(\-x)$, ..., $\psi_r(\-x)$ in $\mathcal V$,
which can be used to represent approximate solutions. 
For convenience's sake, we define the following vector-valued function (feature map):
$$\boldsymbol{\psi} = \begin{bmatrix} \psi_1 \\ \vdots \\ \psi_r \end{bmatrix} \in \mathcal V^r,$$ 
where $\mathcal V^r$ is the $r$-fold Cartesian product of $\mathcal V$. 
Let $\Phi(\cdot;\bdtheta): \mathcal V^r \to \mathcal V$ denote a (possibly nonlinear) operator parameterized by
$\bdtheta \in \Theta \subset \mathbb R^\ell$.
We define the approximate solution manifold as
\[
\mathcal M_{\Phi} := \{\, \Phi(\bpsi;\bdtheta) : \bdtheta \in \Theta \,\} \subset \mathcal V.
\]
In this formulation, the finite-dimensional parameter $\boldsymbol{\theta}$ serves as the latent variable that parametrizes the solution manifold.
That is, for any $u(\cdot,\bmu) \in \mathcal M$, we can find a $\bdtheta \in \Theta$ such that
$u(\cdot,\bmu) \approx \Phi(\bpsi;\bdtheta)$.

While a general operator $\Phi$ provides significant flexibility, its unconstrained nature often leads to computational intractability and a lack of structural guarantees (e.g., regularity) essential for solving PDEs. To enforce a more structured mapping, we restrict $\Phi$ to be a Nemytskii (superposition) operator.
Namely let $h(\cdot;\bdtheta): \mathbb{R}^r  \to \mathbb{R}$ be a  function 
parameterized by parameter $\bdtheta \in \Theta$. We construct the operator $\Phi$ as the  pullback of $h(\cdot,\bdtheta)$
by $\bpsi$:
$$\Phi(\boldsymbol{\psi}; \boldsymbol{\theta})(\mathbf{x}) :=
\boldsymbol{\psi}^* h(\cdot; \boldsymbol{\theta})=
h(\boldsymbol{\psi}(\mathbf{x}); \boldsymbol{\theta}).$$
The operator $\Phi$ constructed in this manner is called an (autonomous) Nemytskii operator,
and unlike general operators that may involve non-local dependencies, $\Phi$ acts \emph{pointwise}, meaning
that its output function at a given coordinate $\mathbf{x}$ depends solely on the local value of the feature map $\boldsymbol{\psi}(\mathbf{x})$.
In practice, representing the operator $\Phi$ reduces to representing the finite-dimensional function
$h(\cdot;\boldsymbol{\theta})$, which we approximate using a neural network;  we refer to the resulting
$\Phi$ as a Nemytskii neural operator (NNO).

The most significant feature of NNO is that the mathematical properties of the output function are well controlled by the feature map $\bpsi$, which, as will be shown later, is fixed rather than learned in the online stage. This helps prevent the model from introducing undesirable or unphysical behavior, especially in unsupervised learning.
Most notably Nemytskii operator $\Phi$ preserves the regularity of the feature map $\bpsi$ in the sense that the
regularity of its output function is inherited from the regularity of $\bpsi$ and the
smoothness of the underlying nonlinear function, which is stated in the following proposition \cite[Section~12.4]{Leoni2017Sobolev}:
\begin{proposition}[Regularity Inheritance]\label{prop:regularity}
Let $\Omega \subset \mathbb{R}^d$ be open, let $h:\mathbb{R}^r \times \Theta \to \mathbb{R}$, and fix $\theta \in \Theta$.
Define the Nemytskii (superposition) operator
\[
\Phi(\boldsymbol{\psi}; \bdtheta)(\-x) := h(\boldsymbol{\psi}(\-x); \bdtheta).
\]
Let $1 \le p < \infty$.
If $\boldsymbol{\psi} \in [W^{1,p}(\Omega) \cap L^\infty(\Omega)]^r$ and
$h(\cdot; \bdtheta)$ is locally Lipschitz on $\mathbb{R}^r$, then
\[
\Phi(\boldsymbol{\psi}; \bdtheta) \in W^{1,p}(\Omega) \cap L^\infty(\Omega).
\]
\end{proposition}



Another important property of NNO is its pointwise nature, which implies that both spatial and parametric derivatives admit explicit pointwise representations, as stated in the following proposition.
\begin{proposition}[Pointwise Structure]
Let $\Omega \subset \mathbb{R}^d$ be a bounded domain and $1 \le p <\infty$. 
Let each $\psi_i\in W^{1,p}(\Omega)\cap L^\infty(\Omega)$ for all
$i=1,..,r$.
Let
\[
h : \mathbb{R}^r \times \Theta \to \mathbb{R},
\qquad \Theta \subset \mathbb{R}^\ell,
\]
be locally Lipschitz in the feature variable $\-z\in \mathbb R^{r}$ and differentiable with respect to the
parameter $\boldsymbol{\theta}$ for almost every $\-z\in \mathbb R^r$. 
Define
\[
u_{\bdtheta}(\-x)
=
h(\boldsymbol{\psi}(\-x);\boldsymbol{\theta}),
\qquad \-x \in \Omega .
\]
Define the (weak) Jacobian matrix of $\boldsymbol{\psi}$ by
\[
J_{\boldsymbol{\psi}}(\-x)
=
\begin{bmatrix}
\nabla \psi_1(\-x)^{\top} \\
\vdots \\
\nabla \psi_r(\-x)^{\top}
\end{bmatrix}
\in \mathbb{R}^{r \times d}.
\]
Then the following hold.
\begin{enumerate}
\item 
The weak spatial gradient of $u_{\boldsymbol{\theta}}$ is given by
\[
\nabla u_{\boldsymbol{\theta}}(\-x)
=
J_{\boldsymbol{\psi}}(\-x)^{\top}
\nabla_z h\big(\boldsymbol{\psi}(\-x);\boldsymbol{\theta}\big),
\qquad \text{a.e. } \-x \in \Omega .
\]

\item 
For each $i = 1,\dots,\ell$, the partial derivative of $u_{\boldsymbol{\theta}}$
with respect to $\theta_i$ is given pointwise by
\[
\frac{\partial u_{\boldsymbol{\theta}}(\-x)}{\partial \theta_i}
=
\frac{\partial h}{\partial \theta_i}
\big(\boldsymbol{\psi}(\-x);\boldsymbol{\theta}\big),
\qquad \text{a.e. } \-x \in \Omega .
\]  
\end{enumerate}
\end{proposition}
This result is a direct consequence of the Sobolev chain rule for Nemytskii operators; see \cite[Chapter~II]{appell1990nonlinear}.
In our method, gradients of the approximate solution with respect to both the spatial variable $\-x$ and the parameter $\boldsymbol{\theta}$ are evaluated repeatedly. The explicit pointwise formulas provided above therefore lead to substantial computational savings.

Finally, we provide a remark clarifying how the assumptions in the preceding propositions are satisfied by common NN models.
\begin{rem}
    The assumptions imposed on the nonlinear function $h$ in the preceding propositions are mild,
    and  standard feedforward NN models with commonly used activation functions, such as ReLU, leaky ReLU, sigmoid, tanh, Swish and softplus, satisfy the assumptions. 
    It is interesting to note  that these assumptions are also satisfied by networks with nonsmooth activation functions, such as ReLU and leaky ReLU.
\end{rem}

\subsection{Offline training} \label{sec:offline}

The offline stage of our method is very similar to that in the non-intrusive RB methods.
Recall that in the non-intrusive RB framework, the main task of the offline stage is two-fold: first, it constructs the reduced bases from the snapshots; second, it builds the regression model from the parameters to the RB coefficients, which can then be used to obtain the reduced solution directly.
This task remains the same in the offline stage of our method, but there are two main differences.
First, the dimension reduction process is nonlinear and data-driven.
Second, unlike the standard and non-intrusive RB methods, which assume that the reduced solution is a linear combination of the basis functions, our approach employs a nonlinear mapping from the basis functions to the reduced solution.

The starting point of our method is the same as for the standard reduced basis methods: 
generate a set of  snapshots and obtain RB functions $\{\psi'_1, \ldots, \psi'_{r'}\}$ via POD of the snapshots. 
Here, we refer to $\{\psi'_1, \ldots, \psi'_{r'}\}$ as the ``initial feature functions'', 
and the actual feature functions, denoted by $\{\psi_1,\ldots,\psi_{r}\}$, are obtained through a nonlinear mapping of the initial ones:
\begin{equation} \label{e:map:Phi}
[\psi_1,\ldots,\psi_r]= \Psi(\psi'_1,\ldots,\psi'_{r'}),
\end{equation}
where $r \leq r'$ is the number of the actual feature functions. 
The purpose of mapping $\Psi$ is dimension reduction (DR) or feature extraction, and here we  restrict it to be a NNO:
\begin{equation} \label{e:Phi_NNO}
\bpsi(\-x) = \Psi(\bpsi'(\-x);\bdtheta_\mathrm{DR})(\-x) = \-g(\bpsi'(\-x);\bdtheta_\mathrm{DR}),
\end{equation}
where 
$\boldsymbol{\psi}'(\-x) = [\psi'_1(\-x),\dots,\psi'_{r'}(\-x)]^{\top} \in \mathcal V^{r'}$,
and $\-g: \mathbb R^{r'}\rightarrow \mathbb R^r$ is represented by an NN model with parameter $\boldsymbol{\theta}_\mathrm{DR}$. 

\begin{figure}[h]
    \centering
    \includegraphics[width=.8\textwidth]{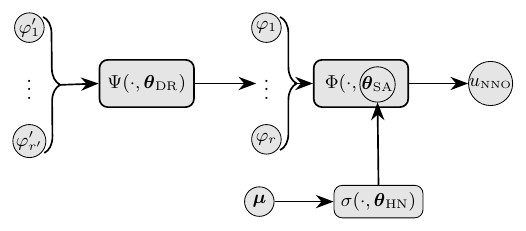}
    \caption{The offline training procedure and the network architecture.}
    \label{fig:offline}
\end{figure}

Next, we assume that solutions can be represented by a NNO of the actual  feature functions:
\begin{equation}
u(\-x) \approx \Phi(\bpsi(\-x);\bdtheta_\mathrm{SA})(\-x)
= h(\bpsi(\-x);\bdtheta_\mathrm{SA}). \label{e:sol_rec}
\end{equation}
where $h(\cdot;\bdtheta_\mathrm{SA}): \mathbb R^r \rightarrow \mathbb R$
is  an NN model parametrized by $\boldsymbol{\theta}_\mathrm{SA}$. 
Network $h$ should be kept shallow and narrow to maintain simplicity and efficiency.
Recall that the solution depends on the PDE parameter $\boldsymbol{\mu}$,
and in our method, parameter $\boldsymbol{\mu}$ enters the formulation through the network parameter
$\boldsymbol{\theta}_\mathrm{SA}$. Specifically we construct a hypernet (HN) of $h(\cdot;\bdtheta_\mathrm{SA})$
by writing 
$\boldsymbol{\theta}_\mathrm{SA}$ as:
\begin{equation}
\boldsymbol{\theta}_\mathrm{SA}= \sigma(\boldsymbol{\mu};\boldsymbol{\theta}_\mathrm{HN}),
\label{e:network:Theta}
\end{equation}
where $\sigma(\cdot;\boldsymbol{\theta}_\mathrm{HN})$ is a neural network
with parameter $\boldsymbol{\theta}_\mathrm{HN}$. 
We emphasize that NN $\-g(\cdot;\bdtheta_\mathrm{DR})$ (and hence operator $\Psi$) does not depend on $\bdmu$, implying that the feature functions $\psi_1$,\ldots, $\psi_r$ are independent of 
$\bdmu$. 
The entire model architecture is illustrated by Figure~\ref{fig:offline}.

To summarize, our model consists of  three NNs: 
the \emph{dimension reduction} model $\Psi$, the \emph{solution approximation} model $\Phi$ and the \emph{hypernet}  model $\sigma$.
For conciseness we refer to the entire NNO model in the offline stage 
as,  
\begin{multline}
u_{\rm off}(\-x,\boldsymbol{\mu};\boldsymbol{\theta}_\mathrm{DR},\bdtheta_\mathrm{HN}) 
= \Phi(\Psi(\bpsi';\bdtheta_\mathrm{DR});\sigma(\bmu;\bdtheta_\mathrm{HN}))(\-x)\\
= h(\-g(\bpsi'(\-x);\bdtheta_\mathrm{DR});\sigma(\bmu;\bdtheta_\mathrm{HN})),
\end{multline}
 as $\bdtheta_\mathrm{SA}$ is the output of hypernet $\sigma$.
Similar to the non-intrusive RB, the network is trained in a supervised manner. 
That is, we generate a new set of $N$ labeled data: $$(\boldsymbol{\mu}^{(1)},u_h(\-x,\bmu^{(1)})),\ldots,(\boldsymbol{\mu}^{(N)},u_h(\-x,\bmu^{(N)})).$$

For training, we use a loss function which encodes
information  about the solution and its first-order derivative. 
We assume that the numerical solutions live in the Sobolev space $W^{1,2}({\Omega})$
and define the offline loss function using the weighted $W^{1,2}(\Omega)$-norm:
\begin{multline}
L_{\rm off}(\boldsymbol{\theta}) = \sum_{i=1}^{N}\left \| u_{\rm off}(\mathbf{x},\boldsymbol{\mu}^{(i)}; 
\boldsymbol{\theta}_\mathrm{DR},\bdtheta_\mathrm{HN}) -  u_h(\mathbf{x},\bmu^{(i)}) \right\|_{L^2(\Omega)}^2 
\\+\lambda_{\nabla}\sum_{i=1}^{N}\left \| \nabla u_{\rm off}(\mathbf{x},\boldsymbol{\mu}^{(i)}; 
\boldsymbol{\theta}_\mathrm{DR},\bdtheta_\mathrm{HN}) -  \nabla u_h(\mathbf{x},\boldsymbol{\mu}^{(i)}) \right\|_{L^2(\Omega)}^2,
\label{e:sobolev}
\end{multline}
where \( \lambda_{\nabla} \) is a weight parameter.
As such, in the offline stage, the NNO model is trained by minimizing the loss function \eqref{e:sobolev}.

\subsection{Online adaptation}
Through the offline stage, we have obtained a model $F$ which can be used directly to predict the solution for a given parameter value $\bdmu^*$. 
However, this model may not be sufficiently accurate,  especially for parameter values in regions with sparse training data.
To further improve the prediction accuracy, we propose to conduct unsupervised learning to adapt the model during the online stage, using a formulation similar to that of PINNs. 

\begin{figure}[h]
    \centering
    \includegraphics[width=.5\textwidth]{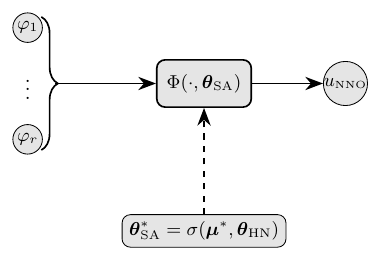}
    \caption{The online training procedure and the network architecture.}
    \label{fig:online}
\end{figure}

In the offline stage, we have trained three networks: the DR network $\Psi$, the SA network
$\Phi$, and the hypernet $\sigma$. 
In the online stage, the DR network $\Psi$ is frozen, which means that the feature functions obtained in the offline stage are used directly.
We seek to adapt the solution reconstruction model $\Phi(\bpsi;\bdtheta_\mathrm{SA})$
in an unsupervised manner. With the feature map $\bpsi$ being fixed, the NNO solution in the online stage can be written as
\begin{equation}
u_\mathrm{on}(\-x;\bdtheta_\mathrm{SA}) = \Phi(\bpsi;\bdtheta_\mathrm{SA})(\-x)
= h(\bpsi(\-x);\bdtheta_\mathrm{SA}), \label{e:sol_online}
\end{equation}
and is trained using the  PINNs loss~\eqref{e:loss_pinn} with a collection of pre-selected test functions.
To accelerate the computation, the optimization problem is initialized 
using the output of the hypernet model obtained at the offline stage:
$$\bdtheta_\mathrm{SA}^*=\sigma(\bdmu^*;\bdtheta_\mathrm{HN}).$$
The online training procedure and the network architecture are shown in Figure~\ref{fig:online}.

It is important to note that online training typically can be carried out very efficiently, owing to several factors.
First, as mentioned in Section \ref{sec:offline}, the network $\Phi$ has been deliberately kept small in both depth and width, 
and so the dimensionality of the parameter space is modest. For example, in Section \ref{sec:burgers}, a network with a single hidden layer of 5 neurons is used to reconstruct the solution of the Burgers' equation.
Second, the evaluation and optimization of the PINN loss requires the computation of the reconstructed NN-based solution and its derivatives (with respect to 
both spatial variables and network parameters) at a large number of collocation points. Following Proposition 2, the NNO formulation provides pointwise explicit expressions for both the spatial and the parameter gradients, considerably reducing the computational cost of optimizing the loss. 
Finally, in most cases, $\bdtheta_\mathrm{SA}^*$ provides a good initial guess for the training procedure, so that the number of iterations 
can be significantly reduced.

\subsection{Connections to PINNs}
We briefly discuss the motivation of the NNO method from the PINN perspective. 
PINNs are often computationally intensive when applied to complex PDE models: firstly, in order to accurately approximate solutions to such models, the neural network typically needs to be large in both depth and width,
and secondly, a substantial number of training iterations is often required to achieve a satisfactory solution of the associated minimization problem.
In what follows, we will discuss how the proposed method can mitigate this computational burden by introducing an offline pretraining stage, which reduces the complexity of the online PINN training.

The goal of PINNs is to construct a scalar-valued function of the physical variables that approximates the true solution.
In this sense, PINNs can be seen as analogous to classical regression problems, where one approximates an underlying target function within a prescribed model class.
In regression, feature maps are often introduced to simplify the approximation of complicated and highly nonlinear target functions by representing them in a more expressive or structured feature space.
In this regard, the main idea of the NNO method is to construct an explicit feature map and to approximate the solution as a scalar-valued function defined on the resulting feature space, rather than directly on the physical space.
The underlying motivation is that, in the feature space, the input–output mapping may admit a simpler representation than in the original physical variables, allowing the solution to be represented using a much smaller network than in standard PINNs.
In our method, such effective feature representations are learned from high-fidelity solution snapshots during the offline stage. During the online stage, the solution approximation is then carried out in the resulting feature space using the precomputed feature representation. In summary, 
an alternative interpretation of the proposed NNO approach is a feature-based PINN approach implemented within an offline–online computational framework.

\section{Numerical experiments}\label{sec:examples}
\numberwithin{equation}{section}
We test the proposed NNO method
in two numerical examples, with implementation details provided in the Appendix \ref{sec:implementation}.

\subsection{Burgers' equation} \label{sec:burgers}
\begin{figure}[h]
    \centering
    \includegraphics[width=.85\textwidth]{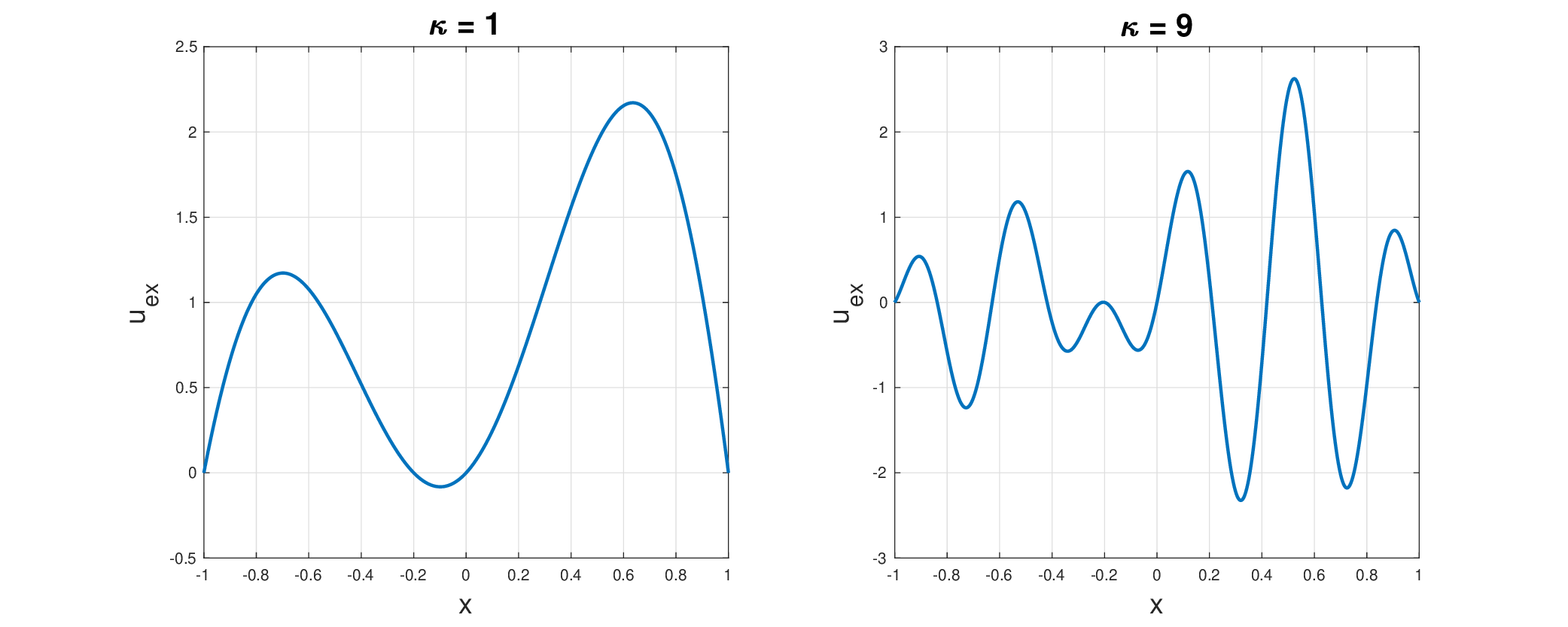}
    \caption{
   Exact solutions of Eq.~\eqref{Ex1 eq} for $\kappa = 1$ (left) and $\kappa = 9$ (right).
}
    \label{fig:burger_sol}
\end{figure}

Our first example is the one-dimensional Burgers' equation over  $\Omega = [-1,1]$:
\begin{equation}\label{Ex1 eq}
\begin{split}
    &u \frac{\partial u}{\partial x} - \frac{\partial^2 u}{\partial x^2}
    = s(x,\boldsymbol{\mu}),  \\
    &u(-1) =   u(1) = 0,
\end{split}
\end{equation}
where $s(x,\bdmu)$ is the source term and $\bdmu$ is a two-dimensional parameter $\boldsymbol{\mu} = [\mu_1,\mu_2]^T$ following a uniform distribution defined on $[1,10]\times[1,10]$.
 In this example we choose the source term to be
\begin{align*}
s(x,\bdmu)&=
(1+\mu_1 x)^2\Bigl[\,6x^2-2
-\frac{2\,\kappa^2\mu_2^2}{9}\,(x^2-1)^2\Bigr]
\cos\!\Bigl(\tfrac{2\kappa\mu_2x}{3}\Bigr)\\[6pt]
&\quad
+\;\frac{2\kappa\mu_2}{3}\,(1+\mu_1 x)\,(x^2-1)\,
\bigl[\mu_1(x^2-1)+2x(1+\mu_1 x)\bigr]\,
\sin\!\Bigl(\tfrac{2\kappa\mu_2x}{3}\Bigr)\\[6pt]
&\quad
-\;(6x^2-2)\,(1+\mu_1 x)^2\,,
\end{align*} 
and it follows that Eq.~\eqref{Ex1 eq} admits the exact solution given by
\begin{equation}\label{1d exact}
   u_{\rm ex}(x,\boldsymbol{\mu}) = (1+\mu_1 x)\sin{(-\kappa\frac{\mu_2 x}{3})}(x^2-1). 
\end{equation}
Note that in Eq.~\eqref{1d exact} $\kappa$ is a positive parameter controlling the ``frequency'' of the solution, i.e., how rapidly it varies with respect to $x$.   
To demonstrate the effect of $\kappa$,  Fig. \ref{fig:burger_sol} shows two sample solutions corresponding to $\kappa = 1$ and $\kappa = 9$ (both with $\mu_1=\mu_2=5$), respectively.
As shown in the figure, the solution exhibits more rapid variation for larger $\kappa$, making it intuitively more challenging to approximate.
In this example, we will evaluate the performance of the proposed method for both $\kappa = 1$ and $\kappa = 9$.

\begin{figure}[h]
    \centering
\centerline {   \includegraphics[width=.5\textwidth]{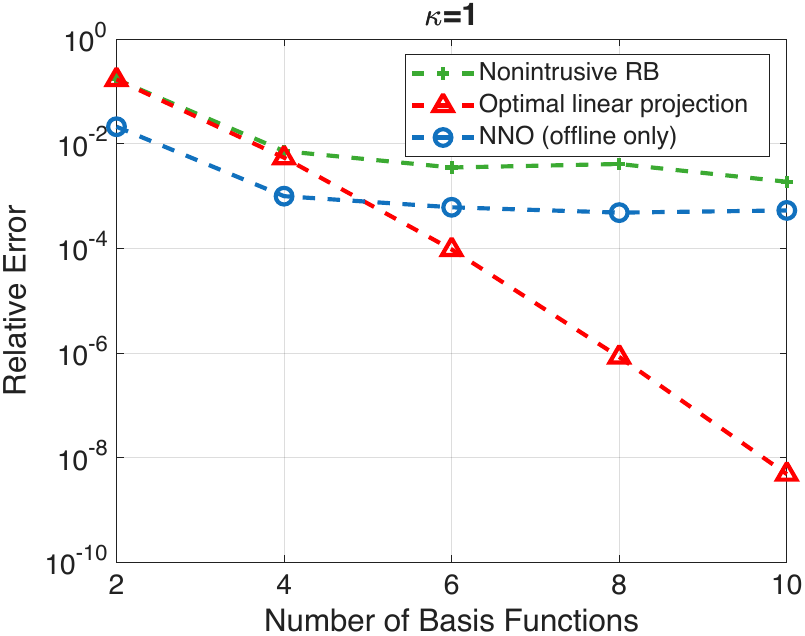}
\includegraphics[width=.5\textwidth]{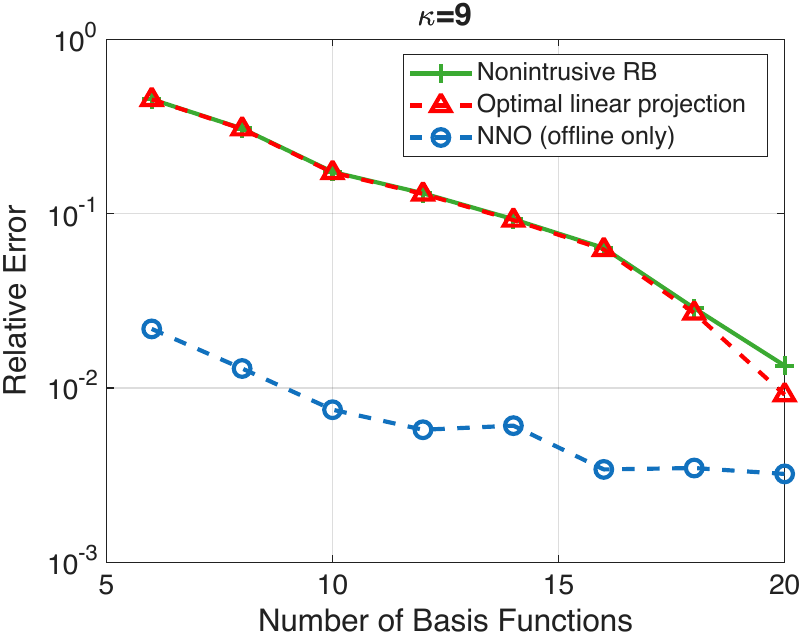}}
    \caption{  Relative errors of the offline-computed solutions  as a function of the number of basis/feature functions. Left for the case $\kappa = 1$, and right for the case $\kappa = 9$.}
    \label{fig:offline_burger}
\end{figure}

In all numerical experiments, we ensure that the set of snapshots is sufficiently large so that it does not become a major limiting factor in performance: the sample size is 100 for $\kappa = 1$, and 1600 for $\kappa = 9$.  
We first examine the NNO model obtained solely from the offline computation.  
In Fig.~\ref{fig:offline_burger}, we plot the relative mean-squared errors of the reduced solutions versus the number of basis functions for both $\kappa = 1$ and $\kappa= 9$. 
The errors are computed using 100 randomly generated test samples. 
In addition to the results of the proposed method, we show two reference results for comparison.  
 The first is the optimal linear projection, obtained by projecting the exact solution onto the reduced basis; it represents the best-case scenario when a linear mapping is used.
The second is the result of the non-intrusive RB method, which, as discussed earlier, is also based on linear combination. 
For the less oscillatory case ($\kappa = 1$), we observe that the error in the optimal linear projection decays rapidly, becoming smaller than that of our method when six or more basis functions are used. This suggests that in such a simple setting, a linear combination is capable of approximating the solution quite well.
Obviously the optimal projection serves as an ideal benchmark, rather than a practical method.  
When comparing the results of the NNO method with those of the linear-combination-based non-intrusive RB, we can see that 
the former consistently yields smaller errors, demonstrating the advantage of incorporating nonlinearity in the reduced model.
The advantage of the nonlinear NNO becomes more evident in the case of $\kappa = 9$.  
In this more challenging setting, our method not only significantly outperforms the non-intrusive RB method, but also achieves considerably lower errors than the optimal linear projection when using up to 20 basis/feature functions.  
This highlights the effectiveness of NNO in capturing complex solution behavior.

\begin{figure}[h]

\centerline{
\includegraphics[width=0.49\linewidth]{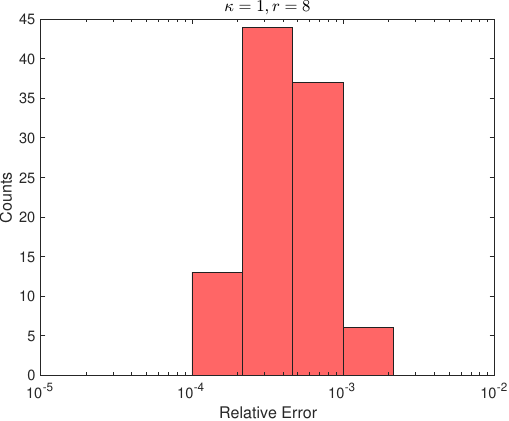} 
\includegraphics[width=0.49\linewidth]{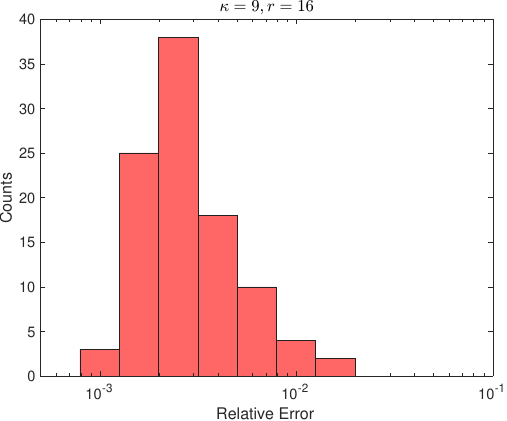} }
\caption{Histograms of the relative errors of the offline-trained model over 100 test samples. Left for the case $\kappa = 1$ using 8 feature functions; right for the case $\kappa= 9$ using 16 feature functions.}
\label{fig:hist}
\end{figure}

\begin{figure}[h]
	\begin{subfigure}{0.49\textwidth}
  \includegraphics[width=1\textwidth]{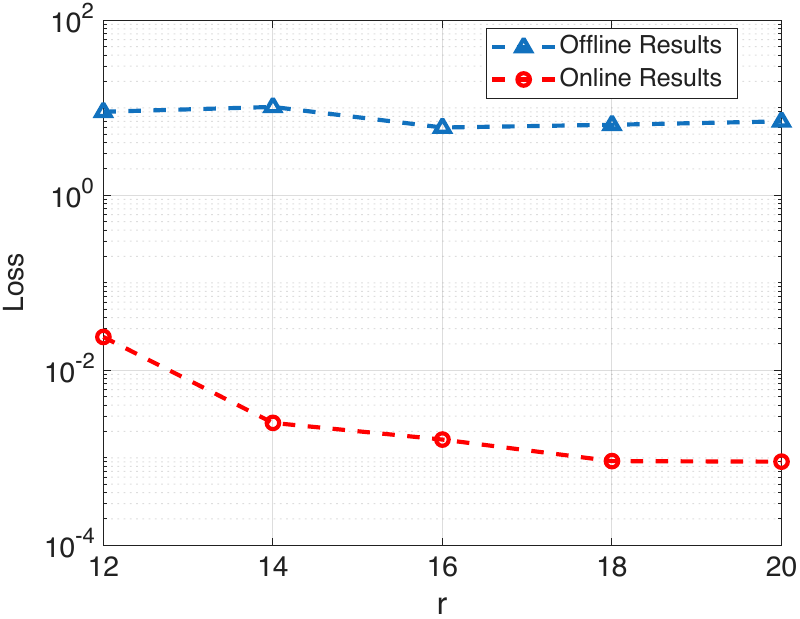}	
	\caption{Online vs Offline: residual losses averaged over all test samples}
\end{subfigure}
\begin{subfigure}{0.49\textwidth}
	\includegraphics[width=1\textwidth]{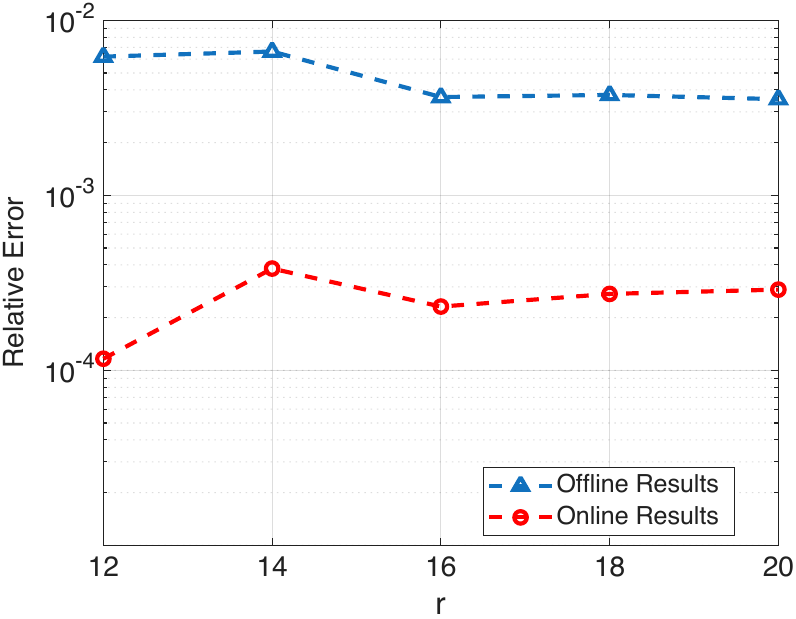}
		\caption{Online vs Offline: relative errors averaged over all test samples}
	\end{subfigure}
	\caption{Relative errors and residual losses of the solutions reconstructed using the offline stage only versus those with online adaptation.}
	\label{fig:on_off}
\end{figure}

Next, we consider the effect of online network adaptation. 
The first question we aim to address is whether such online computation is necessary. 
To this end, we plot in Fig.~\ref{fig:hist} the histogram of the relative errors of the reduced solutions obtained by 
the offline training only.
For $\kappa = 1$, we use $r = 8$ basis/feature functions, and for $\kappa = 9$, we set $r = 16$.
The results show that the offline model performs very well in the simpler case ($\kappa = 1$), but less effectively in the more challenging case ($\kappa = 9$).
If we use $10^{-3}$ as a threshold for acceptable accuracy,  nearly all test cases—94 out of 100—satisfy this criterion for $\kappa = 1$, while for $\kappa = 9$, only 3 out of 100 cases fall below it.
From this we can conclude that online adaptation is particularly important for more complex solutions. 
Therefore,  let us now focus on examining the effect of online adaptation in the $\kappa = 9$ case.  
In Fig.~\ref{fig:on_off} we compare the results of the NNO solutions obtained by the offline training only and those updated with online adaptation.  Note that in the online adaptation stage, the total number of training iterations is fixed to 100.
The left panel shows the average residual loss over 100 test samples before and after online adaptation.
The right panel compares the average relative errors of the same set of samples before and after online adaptation.
As expected, the left panel shows that online training substantially reduces the residual loss, demonstrating the effectiveness of residual loss–driven unsupervised learning.
More importantly, the right panel shows that this reduction translates into a significant improvement in solution accuracy, by approximately one order of magnitude.
These results demonstrate that online adaptation is highly effective in improving the accuracy of the NNO solution, particularly for parameter values that are insufficiently represented in the offline training data.

\subsection{Nonlinear heat conduction problem}\label{lid nobc sec}

We consider a 3-dimensional steady--state nonlinear heat conduction problem
\begin{equation}
\label{eq:example_gov}
\nabla \cdot \bigl(k(u)\nabla u\bigr)
+
S(\mathbf{x};\boldsymbol{\mu})
=
0,
\qquad
\mathbf{x} = (x_1,x_2,x_3) \in \Omega,
\end{equation}
posed on the cubic domain
\[
\Omega = (0,L)^3 \subset \mathbb{R}^3,
\]
with homogeneous Dirichlet boundary conditions
\begin{equation}
\label{eq:example_bc}
u = 0
\qquad
\text{on } \partial \Omega.
\end{equation}
Here, $u(\mathbf{x})$ denotes the temperature field,
and the thermal conductivity is temperature-dependent and given by
\begin{equation}\label{eq:k}
k(u) = k_0(1 + \beta u),
\end{equation}
where $k_0>0$ and $\beta \ge 0$ are constants.
Here $k_0 > 0$ is the reference conductivity and
$\beta \ge 0$ controls the strength of the nonlinearity.
When $\beta = 0$, the problem reduces to a linear Poisson equation,
whereas $\beta > 0$ introduces a nonlinear diffusion effect.
The volumetric heat source is spatially varying but independent of the
temperature:
\begin{equation}\label{eq:source}
S(\mathbf{x};\boldsymbol{\mu})
=
\sum_{i=1}^{N_S} \mu_i
\sin\!\left(\frac{i\pi x_1}{L}\right)
\sin\!\left(\frac{i\pi x_2}{L}\right)
\sin\!\left(\frac{i\pi x_3}{L}\right),
\end{equation}
where
\[
\boldsymbol{\mu} = (\mu_1, \dots, \mu_{N_S}),
\]
is a vector of random parameters. In this example, we take $k_0=1$, $\beta=2$, $L=2$, $N_S=12$ and 
$\boldsymbol{\mu}$ to follow a uniform distribution $U[-1,1]^{N_S}$.
High-fidelity FEM snapshots are obtained using eight-node linear hexahedral elements and a Newton-Raphson solver, with the numerical configuration tuned to ensure the relative error is controlled below $10^{-6}$. 
We generate 1000 snapshots for the offline training, for both the basis generation and network training.

\begin{figure} 
        \centerline{
        \includegraphics[width=.43\linewidth]{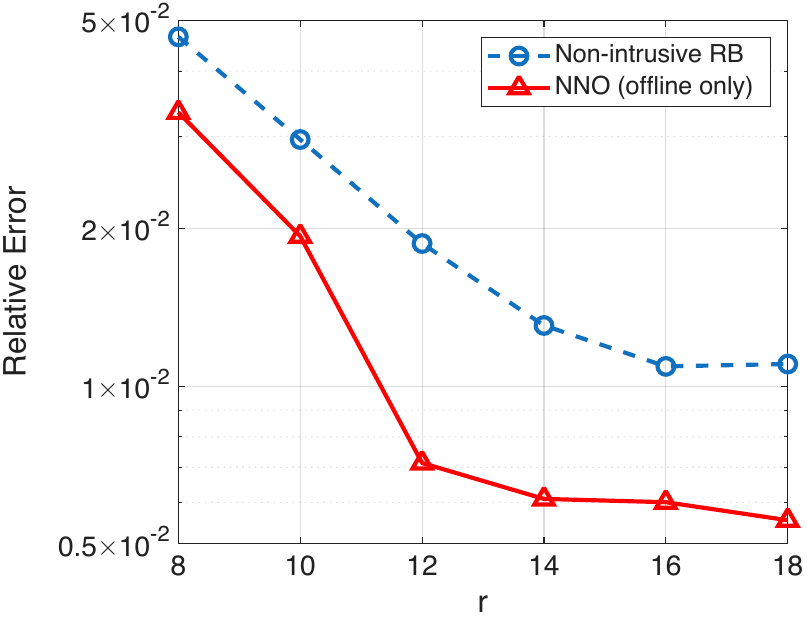}
 \includegraphics[width=.43\linewidth]{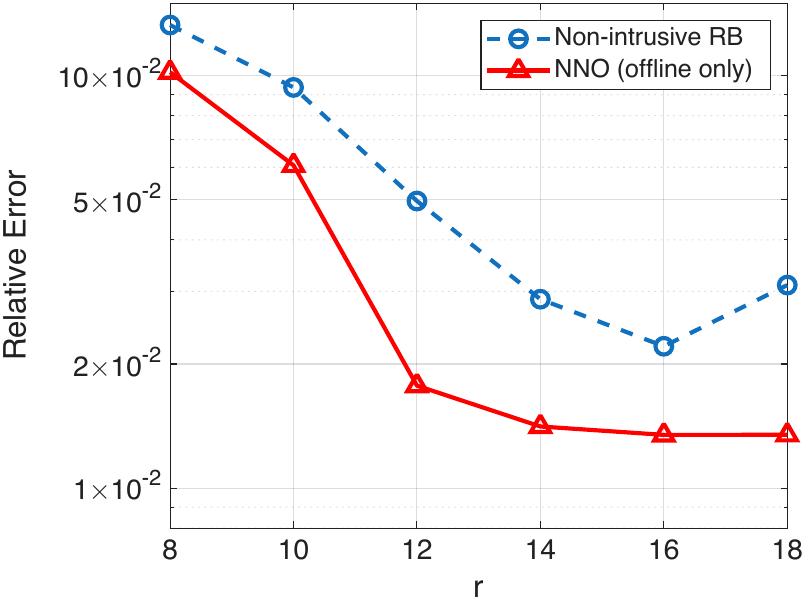}
    }
    \caption{Relative errors of the offline-trained model as a function of the number of basis/feature functions. 
    Left: Errors averaged across all 100 test samples; Right: Errors averaged across the 10 worst cases.}
    \label{fig:heat_offline}
\end{figure}
\begin{figure} 
    \centerline{\includegraphics[width=.45\linewidth]{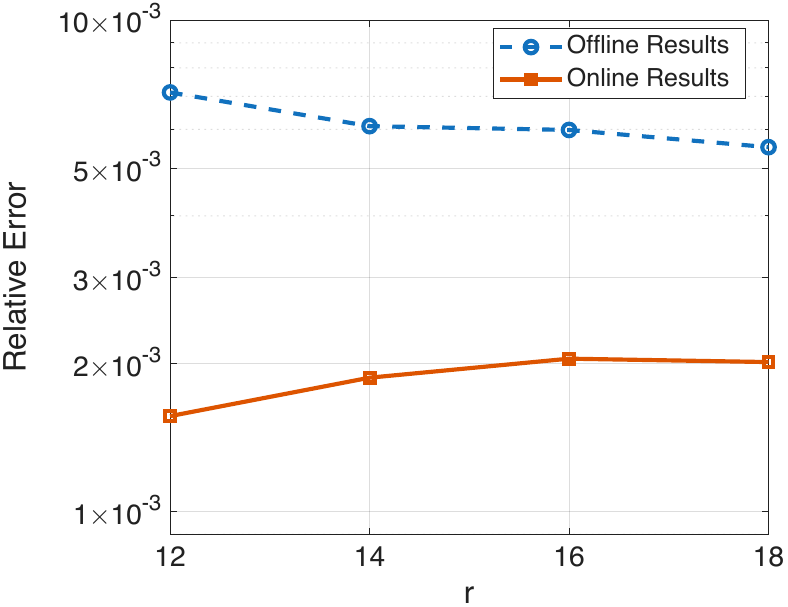}
\includegraphics[width=.45\linewidth]{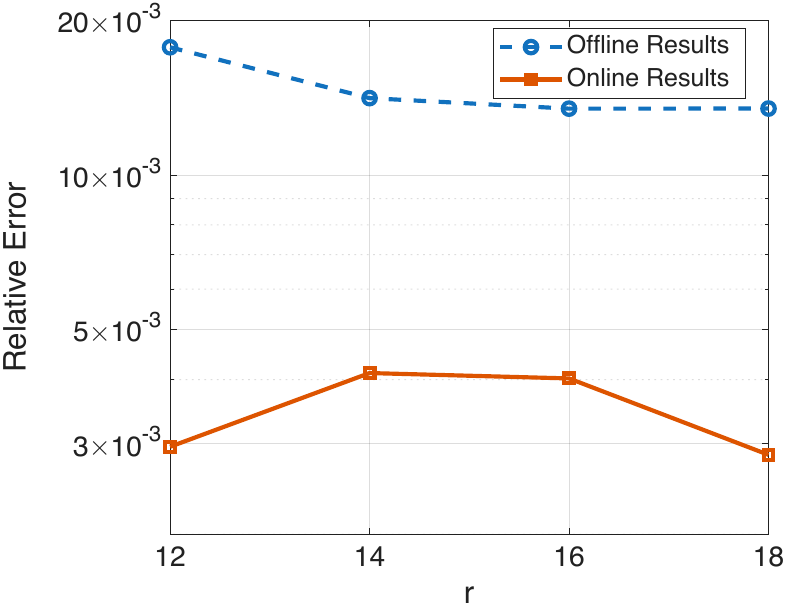}}
    \caption{Relative errors of the solutions reconstructed using the offline stage only versus those with online adaptation.
    Left: Errors averaged across all 100 test samples; Right: Errors averaged across the 10 worst cases.}
    \label{fig:heat_offline_online}
\end{figure}

We first consider the offline NNO model and compare its results with the non-intrusive RB method, using 100 randomly generated realizations of the PDE parameters.
In Fig.~\ref{fig:heat_offline} (left), we present the relative errors of the NNO method and the non-intrusive RB method as functions of the number of basis functions, averaged over 100 randomly generated test samples. The results show that the nonlinear RB method consistently outperforms the non-intrusive RB approach across all tested basis sizes, indicating that the nonlinear construction is more effective at capturing the intrinsic structure of the solution manifold in this example.
Figure~\ref{fig:heat_offline} (right) further reports the worst 10 cases among the 100 test samples. Although the overall trends remain similar, the errors in these challenging instances are substantially larger than the averaged results. In particular, even for relatively large basis dimensions, the worst-case error of NNO remains on the order of $10^{-2}$, while that of the non-intrusive RB method is considerably higher. This observation indicates that the offline-constructed 
NN models are not sufficiently robust to uniformly approximate all parameter realizations.


Next, we examine the results of the online adaptation using the same 100 test samples that were used to evaluate the offline results.
In Fig.~\ref{fig:heat_offline_online}, we plot the relative errors of the offline-only predictions
and those with online adaptation,
against the number of  RB functions. The results are shown for  the average error over all test samples (left plot) and the average error over the 10 worst-case samples having
the largest offline prediction errors (right plot). 
As shown in the figures, the online adaptation substantially improves the accuracy of the NNO models.
The improvement is particularly significant for the worst-case samples, where the offline errors are largest. In these cases, the average error is reduced from above $1.2\times10^{-2}$ to less than $4\times10^{-3}$, bringing it well below the $10^{-2}$ level.

\begin{figure}
    \centerline{
    \includegraphics[width=0.85\linewidth]{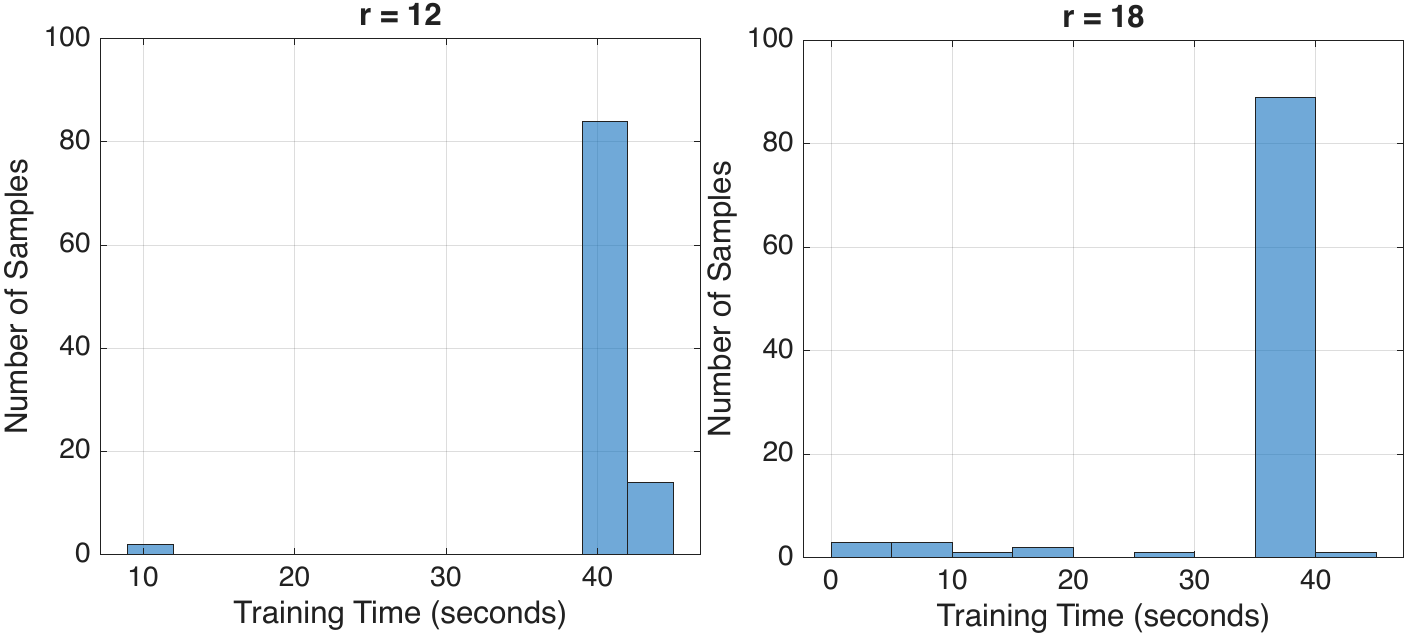}}
    \caption{Histograms of online training times for 100 test examples. Left: the results for 
     12 feature functions.
    Right: the results for 
    18 feature functions.}
    \label{fig:hist_online_time}
\end{figure}

Finally, we discuss the computational cost of solving the parametrized PDE problem in this example. 
We report wall-clock times based on our specific network architecture, as detailed in Appendix \ref{sec:implementation}.
Generating a single high-fidelity snapshot typically takes approximately two hours. 
The offline training phase requires about 20 hours for configurations with both 12 and 18 feature functions.
The time for online adaptation varies across different test samples,
and Fig.~\ref{fig:hist_online_time} shows the histograms of online training times for the 100 test examples, for two 
different numbers of feature functions 12 and 18.
The figure shows that, for both configurations, the online training takes less than 50 
 seconds, which is considerably more efficient than using the conventional solver. 
We note that, although offline training takes significantly longer than online training, the primary offline computational cost comes from high-fidelity snapshot generation, similar to non-intrusive RB methods.

\section{Conclusions} \label{sec:conclusion}

In this work, we have introduced an NNO framework as a nonlinear generalization of RB methods for parametrized steady-state PDEs. The method replaces the standard linear combination of basis functions with a structured, nonlinear mapping realized by a pointwise Nemytskii operator acting on fixed feature functions. By integrating nonlinear dimension reduction in the offline stage with lightweight, physics-informed online adaptation, the NNO framework accurately captures complex solution manifolds that cannot be well approximated
by linear subspaces. Our numerical experiments on the one-dimensional Burgers' equation and a three-dimensional nonlinear heat conduction problem demonstrate that the proposed method consistently outperforms the linear RB approach. 
Specifically, the online adaptation stage proves highly effective, yielding substantial improvements in prediction accuracy, particularly in parameter regimes where training data are limited.

There are several possible directions for extending and improving the proposed method. {First and foremost}, while in this work we have restricted ourselves to steady-state PDE problems, 
a natural extension is to consider unsteady, or time-dependent, PDE systems, where the reduced model must evolve in time. This could be achieved by incorporating
several techniques that have been developed for time-dependent PDEs or dynamical systems  \cite{lee2020model,regazzoni2019machine}. 
 {Second}, the method can be linked to recent developments in {operator-based neural networks}, such as DeepONets \cite{lu2019deeponet,lu2021learning} and neural operators \cite{kovachki2023neural}, which may offer improved generalization capabilities. {Third}, in the present formulation, the dimension reduction network is 
 chosen to be independent of the model parameters,
 and so are the obtained basis functions. 
 Alternatively, introducing {parameter-dependent reduced bases or nonlinear dimension reduction} 
 could enable the model to better capture local features of the solution manifold and handle problems with strong transport or discontinuities \cite{nair2019transported,rim2023manifold,chen2024tgpt}. {Finally}, the proposed approach is well suited for such
 real-world {applications as Bayesian inference or inverse problems}, where its non-intrusive and adaptive nature may significantly accelerate posterior sampling or optimization.  These directions will further broaden the applicability and impact of the proposed method across scientific and engineering domains, and we will explore them in future work.

\appendix

\section{Implementation details} \label{sec:implementation}

In this section, we briefly provide the implementation details of our numerical experiments.
All experiments were conducted on a computing node equipped with dual-socket Intel(R) Xeon(R) Platinum 8360Y CPUs (2.40 GHz, 72 cores in total), with 28 cores allocated for the computations. The system also includes an NVIDIA A100-SXM4-40GB GPU and 503 GiB of RAM.

In both numerical examples, we use  the same network architectures, which are detailed below. 
The \emph{dimension reduction} model $\-g$ in~\eqref{e:map:Phi} consists of one fully connected layer with Swish activation \cite{nwankpa2018activation}, followed by an output layer with no activation. This network has $r'$ input neurons, $r$ neurons at the hidden layer and $r$ outputs. We fix $r'=20$ in our experiments and study the performance of the model for different values of $r$. 
 The \emph{solution reconstruction} model $h$ (see~\eqref{e:sol_rec}) is a one hidden layer network with $r$ inputs, $5$ hidden neurons and one output. The activation used is Tanh~\cite{nwankpa2018activation}. We emphasize that this network is deliberately kept simple as it needs to be trained online. 
 The input dimension of the hypernet model $\sigma$ (see~\eqref{e:network:Theta}) is equal to the dimension of the parameter domain; this model has 
 four hidden layers with $100$ neurons each. The activation after each hidden layer is Swish. The number of outputs of $\sigma$ is the total number of network parameters in $h$. 
 Specifically, $h$ is represented by 
 two weight matrices and two bias vectors, so the total number of network parameters is $5r + 11$.





\begin{singlespace}
	\bibliographystyle{siam}

	\bibliography{ref.bib}

@article {cai2021physics,
    AUTHOR = {Cai, Shengze and Mao, Zhiping and Wang, Zhicheng and Yin,
              Minglang and Karniadakis, George Em},
     TITLE = {Physics-informed neural networks ({PINN}s) for fluid
              mechanics: a review},
   JOURNAL = {Acta Mech. Sin.},
  FJOURNAL = {Acta Mechanica Sinica},
    VOLUME = {37},
      YEAR = {2021},
    NUMBER = {12},
     PAGES = {1727--1738},
      ISSN = {0567-7718,1614-3116},
   MRCLASS = {76-10},
  MRNUMBER = {4412280},
       DOI = {10.1007/s10409-021-01148-1},
       URL = {https://doi.org/10.1007/s10409-021-01148-1},
}

@article {carlberg2011efficient,
    AUTHOR = {Carlberg, Kevin and Bou-Mosleh, Charbel and Farhat, Charbel},
     TITLE = {Efficient non-linear model reduction via a least-squares
              {P}etrov-{G}alerkin projection and compressive tensor
              approximations},
   JOURNAL = {Internat. J. Numer. Methods Engrg.},
  FJOURNAL = {International Journal for Numerical Methods in Engineering},
    VOLUME = {86},
      YEAR = {2011},
    NUMBER = {2},
     PAGES = {155--181},
      ISSN = {0029-5981,1097-0207},
   MRCLASS = {76M25},
  MRNUMBER = {2814379},
       DOI = {10.1002/nme.3050},
       URL = {https://doi.org/10.1002/nme.3050},
}

@article {casenave2015nonintrusive,
    AUTHOR = {Casenave, Fabien and Ern, Alexandre and Leli\`evre, Tony},
     TITLE = {A nonintrusive reduced basis method applied to aeroacoustic
              simulations},
   JOURNAL = {Adv. Comput. Math.},
  FJOURNAL = {Advances in Computational Mathematics},
    VOLUME = {41},
      YEAR = {2015},
    NUMBER = {5},
     PAGES = {961--986},
      ISSN = {1019-7168,1572-9044},
   MRCLASS = {65N30 (65J05 76Q05)},
  MRNUMBER = {3428555},
MRREVIEWER = {C.\ Ilioi},
       DOI = {10.1007/s10444-014-9365-0},
       URL = {https://doi.org/10.1007/s10444-014-9365-0},
}

@article {chaturantabut2010nonlinear,
    AUTHOR = {Chaturantabut, Saifon and Sorensen, Danny C.},
     TITLE = {Nonlinear model reduction via discrete empirical
              interpolation},
   JOURNAL = {SIAM J. Sci. Comput.},
  FJOURNAL = {SIAM Journal on Scientific Computing},
    VOLUME = {32},
      YEAR = {2010},
    NUMBER = {5},
     PAGES = {2737--2764},
      ISSN = {1064-8275,1095-7197},
   MRCLASS = {65L99 (35J60 35K55 65F99 65M99)},
  MRNUMBER = {2684735},
MRREVIEWER = {Sebastian\ Reich},
       DOI = {10.1137/090766498},
       URL = {https://doi.org/10.1137/090766498},
}

@article {chen2021physics,
    AUTHOR = {Chen, Wenqian and Wang, Qian and Hesthaven, Jan S. and Zhang,
              Chuhua},
     TITLE = {Physics-informed machine learning for reduced-order modeling
              of nonlinear problems},
   JOURNAL = {J. Comput. Phys.},
  FJOURNAL = {Journal of Computational Physics},
    VOLUME = {446},
      YEAR = {2021},
     PAGES = {Paper No. 110666, 28},
      ISSN = {0021-9991,1090-2716},
   MRCLASS = {35A25 (65M70 65M99)},
  MRNUMBER = {4308798},
       DOI = {10.1016/j.jcp.2021.110666},
       URL = {https://doi.org/10.1016/j.jcp.2021.110666},
}

@article {chen2024gpt,
    AUTHOR = {Chen, Yanlai and Koohy, Shawn},
     TITLE = {G{PT}-{PINN}: generative pre-trained physics-informed neural
              networks toward non-intrusive meta-learning of parametric
              {PDE}s},
   JOURNAL = {Finite Elem. Anal. Des.},
  FJOURNAL = {Finite Elements in Analysis and Design},
    VOLUME = {228},
      YEAR = {2024},
     PAGES = {Paper No. 104047, 15},
      ISSN = {0168-874X,1872-6925},
   MRCLASS = {65M70 (68T07)},
  MRNUMBER = {4649937},
       DOI = {10.1016/j.finel.2023.104047},
       URL = {https://doi.org/10.1016/j.finel.2023.104047},
}

@article {chen2024tgpt,
    AUTHOR = {Chen, Yanlai and Ji, Yajie and Narayan, Akil and Xu, Zhenli},
     TITLE = {T{GPT}-{PINN}: nonlinear model reduction with transformed
              {GPT}-{PINN}s},
   JOURNAL = {Comput. Methods Appl. Mech. Engrg.},
  FJOURNAL = {Computer Methods in Applied Mechanics and Engineering},
    VOLUME = {430},
      YEAR = {2024},
     PAGES = {Paper No. 117198, 20},
      ISSN = {0045-7825,1879-2138},
   MRCLASS = {65M70 (65M06)},
  MRNUMBER = {4770265},
       DOI = {10.1016/j.cma.2024.117198},
       URL = {https://doi.org/10.1016/j.cma.2024.117198},
}

@article {fresca2021comprehensive,
    AUTHOR = {Fresca, Stefania and Ded{\'e}, Luca and Manzoni, Andrea},
     TITLE = {A comprehensive deep learning-based approach to reduced order
              modeling of nonlinear time-dependent parametrized {PDE}s},
   JOURNAL = {J. Sci. Comput.},
  FJOURNAL = {Journal of Scientific Computing},
    VOLUME = {87},
      YEAR = {2021},
    NUMBER = {2},
     PAGES = {Paper No. 61, 36},
      ISSN = {0885-7474,1573-7691},
   MRCLASS = {65M60 (65M99 68T07)},
  MRNUMBER = {4244919},
MRREVIEWER = {Yinhua\ Xia},
       DOI = {10.1007/s10915-021-01462-7},
       URL = {https://doi.org/10.1007/s10915-021-01462-7},
}

@article {fresca2022pod,
    AUTHOR = {Fresca, Stefania and Manzoni, Andrea},
     TITLE = {P{OD}-{DL}-{ROM}: enhancing deep learning-based reduced order
              models for nonlinear parametrized {PDE}s by proper orthogonal
              decomposition},
   JOURNAL = {Comput. Methods Appl. Mech. Engrg.},
  FJOURNAL = {Computer Methods in Applied Mechanics and Engineering},
    VOLUME = {388},
      YEAR = {2022},
     PAGES = {Paper No. 114181, 27},
      ISSN = {0045-7825,1879-2138},
   MRCLASS = {65M60 (65M99)},
  MRNUMBER = {4327361},
       DOI = {10.1016/j.cma.2021.114181},
       URL = {https://doi.org/10.1016/j.cma.2021.114181},
}

@article {guo2019data,
    AUTHOR = {Guo, Mengwu and Hesthaven, Jan S.},
     TITLE = {Data-driven reduced order modeling for time-dependent
              problems},
   JOURNAL = {Comput. Methods Appl. Mech. Engrg.},
  FJOURNAL = {Computer Methods in Applied Mechanics and Engineering},
    VOLUME = {345},
      YEAR = {2019},
     PAGES = {75--99},
      ISSN = {0045-7825,1879-2138},
   MRCLASS = {62M99 (62F15 65M60)},
  MRNUMBER = {3880138},
       DOI = {10.1016/j.cma.2018.10.029},
       URL = {https://doi.org/10.1016/j.cma.2018.10.029},
}

@article {guo2018reduced,
    AUTHOR = {Guo, Mengwu and Hesthaven, Jan S.},
     TITLE = {Reduced order modeling for nonlinear structural analysis using
              {G}aussian process regression},
   JOURNAL = {Comput. Methods Appl. Mech. Engrg.},
  FJOURNAL = {Computer Methods in Applied Mechanics and Engineering},
    VOLUME = {341},
      YEAR = {2018},
     PAGES = {807--826},
      ISSN = {0045-7825,1879-2138},
   MRCLASS = {65N30 (65N99)},
  MRNUMBER = {3845646},
       DOI = {10.1016/j.cma.2018.07.017},
       URL = {https://doi.org/10.1016/j.cma.2018.07.017},
}

@article {hesthaven2018non,
    AUTHOR = {Hesthaven, J. S. and Ubbiali, S.},
     TITLE = {Non-intrusive reduced order modeling of nonlinear problems
              using neural networks},
   JOURNAL = {J. Comput. Phys.},
  FJOURNAL = {Journal of Computational Physics},
    VOLUME = {363},
      YEAR = {2018},
     PAGES = {55--78},
      ISSN = {0021-9991,1090-2716},
   MRCLASS = {65N30 (65N22)},
  MRNUMBER = {3784416},
       DOI = {10.1016/j.jcp.2018.02.037},
       URL = {https://doi.org/10.1016/j.jcp.2018.02.037},
}

@article {kovachki2023neural,
    AUTHOR = {Kovachki, Nikola and Li, Zongyi and Liu, Burigede and
              Azizzadenesheli, Kamyar and Bhattacharya, Kaushik and Stuart,
              Andrew and Anandkumar, Anima},
     TITLE = {Neural operator: learning maps between function spaces with
              applications to {PDE}s},
   JOURNAL = {J. Mach. Learn. Res.},
  FJOURNAL = {Journal of Machine Learning Research (JMLR)},
    VOLUME = {24},
      YEAR = {2023},
     PAGES = {Paper No. [89], 97},
      ISSN = {1532-4435,1533-7928},
   MRCLASS = {68T07 (35Q30)},
  MRNUMBER = {4582511},
}

@article {lee2020model,
    AUTHOR = {Lee, Kookjin and Carlberg, Kevin T.},
     TITLE = {Model reduction of dynamical systems on nonlinear manifolds
              using deep convolutional autoencoders},
   JOURNAL = {J. Comput. Phys.},
  FJOURNAL = {Journal of Computational Physics},
    VOLUME = {404},
      YEAR = {2020},
     PAGES = {108973, 32},
      ISSN = {0021-9991,1090-2716},
   MRCLASS = {68T05 (65L60)},
  MRNUMBER = {4043884},
       DOI = {10.1016/j.jcp.2019.108973},
       URL = {https://doi.org/10.1016/j.jcp.2019.108973},
}

@article {nair2019transported,
    AUTHOR = {Nair, Nirmal J. and Balajewicz, Maciej},
     TITLE = {Transported snapshot model order reduction approach for
              parametric, steady-state fluid flows containing
              parameter-dependent shocks},
   JOURNAL = {Internat. J. Numer. Methods Engrg.},
  FJOURNAL = {International Journal for Numerical Methods in Engineering},
    VOLUME = {117},
      YEAR = {2019},
    NUMBER = {12},
     PAGES = {1234--1262},
      ISSN = {0029-5981,1097-0207},
   MRCLASS = {65M99 (76B99 76D05)},
  MRNUMBER = {3911941},
       DOI = {10.1002/nme.5998},
       URL = {https://doi.org/10.1002/nme.5998},
}

@article {pichi2024graph,
    AUTHOR = {Pichi, Federico and Moya, Beatriz and Hesthaven, Jan S.},
     TITLE = {A graph convolutional autoencoder approach to model order
              reduction for parametrized {PDE}s},
   JOURNAL = {J. Comput. Phys.},
  FJOURNAL = {Journal of Computational Physics},
    VOLUME = {501},
      YEAR = {2024},
     PAGES = {Paper No. 112762, 24},
      ISSN = {0021-9991,1090-2716},
   MRCLASS = {65M70 (65M99)},
  MRNUMBER = {4689881},
       DOI = {10.1016/j.jcp.2024.112762},
       URL = {https://doi.org/10.1016/j.jcp.2024.112762},
}

@article {raissi2019physics,
    AUTHOR = {Raissi, M. and Perdikaris, P. and Karniadakis, G. E.},
     TITLE = {Physics-informed neural networks: a deep learning framework
              for solving forward and inverse problems involving nonlinear
              partial differential equations},
   JOURNAL = {J. Comput. Phys.},
  FJOURNAL = {Journal of Computational Physics},
    VOLUME = {378},
      YEAR = {2019},
     PAGES = {686--707},
      ISSN = {0021-9991,1090-2716},
   MRCLASS = {65M70 (68T05)},
  MRNUMBER = {3881695},
       DOI = {10.1016/j.jcp.2018.10.045},
       URL = {https://doi.org/10.1016/j.jcp.2018.10.045},
}

@article {regazzoni2019machine,
    AUTHOR = {Regazzoni, F. and Ded\`e, L. and Quarteroni, A.},
     TITLE = {Machine learning for fast and reliable solution of
              time-dependent differential equations},
   JOURNAL = {J. Comput. Phys.},
  FJOURNAL = {Journal of Computational Physics},
    VOLUME = {397},
      YEAR = {2019},
     PAGES = {108852, 26},
      ISSN = {0021-9991,1090-2716},
   MRCLASS = {65P99},
  MRNUMBER = {3990714},
       DOI = {10.1016/j.jcp.2019.07.050},
       URL = {https://doi.org/10.1016/j.jcp.2019.07.050},
}

@article {rim2023manifold,
    AUTHOR = {Rim, Donsub and Peherstorfer, Benjamin and Mandli, Kyle T.},
     TITLE = {Manifold approximations via transported subspaces: model
              reduction for transport-dominated problems},
   JOURNAL = {SIAM J. Sci. Comput.},
  FJOURNAL = {SIAM Journal on Scientific Computing},
    VOLUME = {45},
      YEAR = {2023},
    NUMBER = {1},
     PAGES = {A170--A199},
      ISSN = {1064-8275,1095-7197},
   MRCLASS = {78M34 (35F20 41A46 78M12)},
  MRNUMBER = {4550690},
       DOI = {10.1137/20M1316998},
       URL = {https://doi.org/10.1137/20M1316998},
}

@article {wang2019non,
    AUTHOR = {Wang, Qian and Hesthaven, Jan S. and Ray, Deep},
     TITLE = {Non-intrusive reduced order modeling of unsteady flows using
              artificial neural networks with application to a combustion
              problem},
   JOURNAL = {J. Comput. Phys.},
  FJOURNAL = {Journal of Computational Physics},
    VOLUME = {384},
      YEAR = {2019},
     PAGES = {289--307},
      ISSN = {0021-9991,1090-2716},
   MRCLASS = {80A25 (65M99)},
  MRNUMBER = {3920924},
       DOI = {10.1016/j.jcp.2019.01.031},
       URL = {https://doi.org/10.1016/j.jcp.2019.01.031},
}

@article {xiao2016non,
    AUTHOR = {Xiao, D. and Yang, P. and Fang, F. and Xiang, J. and Pain, C.
              C. and Navon, I. M.},
     TITLE = {Non-intrusive reduced order modelling of fluid-structure
              interactions},
   JOURNAL = {Comput. Methods Appl. Mech. Engrg.},
  FJOURNAL = {Computer Methods in Applied Mechanics and Engineering},
    VOLUME = {303},
      YEAR = {2016},
     PAGES = {35--54},
      ISSN = {0045-7825,1879-2138},
   MRCLASS = {74F10 (76D05)},
  MRNUMBER = {3473456},
       DOI = {10.1016/j.cma.2015.12.029},
       URL = {https://doi.org/10.1016/j.cma.2015.12.029},
}

@article{chen2018greedy,
  title={Greedy nonintrusive reduced order model for fluid dynamics},
  author={W. Chen and J. S. Hesthaven and B. Junqiang and Y. Qiu and Z. Yang and Y. Tihao},
  journal={AIAA Journal},
  volume={56},
  number={12},
  pages={4927--4943},
  year={2018},
  publisher={American Institute of Aeronautics and Astronautics}
}

@book{appell1990nonlinear,
  title={Continuity, Boundedness, and Compactness of Nonlinear Operators},
  author={Appell, J{\"u}rgen and Zabrejko, Petr P.},
  year={1990},
  publisher={Springer-Verlag},
  address={Berlin},
  series={Lecture Notes in Mathematics},
  volume={1400}
}

@book{Leoni2017Sobolev,
  author    = {Giovanni Leoni},
  title     = {A First Course in Sobolev Spaces},
  edition   = {2},
  series    = {Graduate Studies in Mathematics},
  volume    = {181},
  publisher = {American Mathematical Society},
  address   = {Providence, RI},
  year      = {2017},
  isbn      = {978-1-4704-2921-8},
  note      = {See Section 12.4 (Superposition/Nemytskii operators).}
}

@book{goodfellow2016deep,
  title={Deep learning},
  author={I. Goodfellow and Y. Bengio and A. Courville},
  year={2016},
  publisher={MIT press}
}

@book{hesthaven2016certified,
  title={Certified reduced basis methods for parametrized partial differential equations},
  author={J. S. Hesthaven and G. Rozza and B. Stamm},
  volume={590},
  year={2016},
  publisher={Springer}
}

@book{hughes2003finite,
  title={The finite element method: linear static and dynamic finite element analysis},
  author={Hughes, Thomas JR},
  year={2003},
  publisher={Courier Corporation}
}

@article{lu2021learning,
  title={{Learning nonlinear operators via {DeepONet} based on the universal approximation theorem of operators}},
  author={L. Lu and P. Jin and G. Pang and Z. Zhang and G. E. Karniadakis},
  journal={Nature machine intelligence},
  volume={3},
  number={3},
  pages={218--229},
  year={2021},
  publisher={Nature Publishing Group UK London}
}

@article{nwankpa2018activation,
  title={Activation functions: Comparison of trends in practice and research for deep learning},
  author={Nwankpa, C},
  journal={arXiv preprint arXiv:1811.03378},
  year={2018}
}

@article{padula2024brief,
  title={A brief review of reduced order models using intrusive and non-intrusive techniques},
  author={Padula, Guglielmo and Girfoglio, Michele and Rozza, Gianlugi},
  journal={PAMM},
  volume={24},
  number={4},
  pages={e202400210},
  year={2024},
  publisher={Wiley Online Library}
}

@book{quarteroni2015reduced,
  title={Reduced basis methods for partial differential equations: an introduction},
  author={A. Quarteroni and A. Manzoni and F. Negri},
  volume={92},
  year={2015},
  publisher={Springer}
}

@article{yu2019non,
  title={Non-intrusive reduced-order modeling for fluid problems: A brief review},
  author={Yu, Jian and Yan, Chao and Guo, Mengwu},
  journal={Proceedings of the Institution of Mechanical Engineers, Part G: Journal of Aerospace Engineering},
  volume={233},
  number={16},
  pages={5896--5912},
  year={2019},
  publisher={SAGE Publications Sage UK: London, England}
}

@article{lu2019deeponet,
  title={{DeepONet}: Learning nonlinear operators for identifying differential equations based on the universal approximation theorem of operators},
  author={L. Lu and P. Jin and G. E. Karniadakis},
  journal={arXiv preprint arXiv:1910.03193},
  year={2019}
}

@book{demmel1997applied,
  title={Applied numerical linear algebra},
  author={J. W. Demmel},
  year={1997},
  publisher={SIAM}
}
\end{singlespace}

\end{document}